\documentclass[12pt]{article}
\usepackage[utf8]{inputenc}
\usepackage[british]{babel}
\usepackage{cmap}
\usepackage{lmodern}

\usepackage{amssymb, amsmath, amsthm}
\usepackage[a4paper,top=25mm,bottom=25mm,left=25mm,right=25mm]{geometry}
\usepackage{ragged2e}

\usepackage{authblk} % for headings
\usepackage{pifont}
\usepackage{graphicx}
\usepackage[dvipsnames,svgnames,table]{xcolor}
\usepackage[figuresright]{rotating}
\usepackage{xtab} % tackle the long tables
\usepackage{longtable} % tackle the long tables
\usepackage{multirow}
\usepackage{footnote}
\usepackage[stable]{footmisc}
\usepackage{chngpage} % allows for temporary adjustment of side margins
\usepackage{pdflscape} % landscape environment
\usepackage[nottoc,notlot,notlof]{tocbibind} % includes everything in ToC

\usepackage{pgfplots}
\pgfplotsset{every tick label/.append style={font=\footnotesize}}
\pgfplotsset{compat=1.18}
\usepackage{setspace}

\usepackage{array}
\newcolumntype{K}[1]{>{\centering\arraybackslash$}p{#1}<{$}}

\makesavenoteenv{tabular}
\usepackage{tabularx}
\usepackage{booktabs}
\usepackage{threeparttable}
\usepackage[referable]{threeparttablex} % footnotes in tabu
\newcolumntype{R}{>{\raggedleft\arraybackslash}X}
\newcolumntype{L}{>{\raggedright\arraybackslash}X}
\newcolumntype{C}{>{\centering\arraybackslash}X}
\newcolumntype{A}{>{\columncolor{gray!25}}C}
\newcolumntype{a}{>{\columncolor{gray!25}}c}

\newlength{\tablen}

\usepackage{dcolumn} % alignment to decimal points
\newcolumntype{.}{D{.}{.}{-1}}

\usepackage{tikz}
\usetikzlibrary{arrows, calc, matrix, patterns, positioning, trees}
\usepackage[semicolon]{natbib} % numbers,
\usepackage[hyphens]{xurl}
\usepackage{microtype}
\usepackage[justification=centering]{caption} % [justification=centering]

\usepackage[labelformat=simple]{subcaption}

\DeclareCaptionLabelFormat{parenthesis}{(#2)}
\makeatletter
\renewcommand\p@subfigure{\arabic{figure}.}
\makeatother

\DeclareCaptionLabelFormat{parenthesis}{(#2)}
\makeatletter
\renewcommand\p@subtable{\arabic{table}.}
\makeatother

\usepackage{enumitem}

\setlist[itemize]{leftmargin=2.5\parindent}
\setlist[enumerate]{leftmargin=2.5\parindent}

\usepackage{hyperref} % [hidelinks]
\hypersetup{
  colorlinks   = true,    		% Colours links instead of ugly boxes
  urlcolor     = blue,    		% Colour for external hyperlinks
  linkcolor    = blue,    		% Colour of internal links
  citecolor    = ForestGreen	% Colour of citations
}

\def\addlegendimage{\csname pgfplots@addlegendimage\endcsname}

\theoremstyle{plain}
\theoremstyle{definition}
\newtheorem{example}{Example}%[section]

\theoremstyle{remark}

\makeatletter
\let\@fnsymbol\@alph
\makeatother

\def\keywords{\vspace{.5em} % Add keywords
{\noindent \textit{Keywords}: }}

\def\AMS{\vspace{.5em} % Add keywords
{\noindent \textbf{\emph{MSC} class}: }}

\def\JEL{\vspace{.5em} % Add keywords
{\noindent \textbf{\emph{JEL} classification number}: }}

\title{How to optimise tournament draws: \\ The case of the FIFA World Cup}
\author{\href{https://sites.google.com/view/laszlocsato}{L\'aszl\'o Csat\'o}\thanks{~E-mail: \emph{laszlo.csato@sztaki.hun-ren.hu}} }
\affil{Institute for Computer Science and Control (SZTAKI) \\
Hungarian Research Network (HUN-REN) \\
Laboratory on Engineering and Management Intelligence \\
Research Group of Operations Research and Decision Systems}
\affil{Corvinus University of Budapest (BCE) \\
Institute of Operations and Decision Sciences \\
Department of Operations Research and Actuarial Sciences}
\affil{Budapest, Hungary}
\date{\today}

\def\Dedication{
{\noindent
``\emph{In contrast with this, in any other science we consider the object as it has to be represented in accordance with data of intuition; there is a limitless web of intuitions, and therefore of objects of thought, so that the science can never achieve absolute completeness, but can be endlessly extended, as in pure mathematics and the empirical doctrine of Nature.}''\footnote{~Source: \url{https://www.earlymoderntexts.com/assets/pdfs/kant1786.pdf}.}
}
\vspace{0.25cm}

\flushright
\noindent (Immanuel Kant: \emph{Metaphysical Foundations of Natural Science})

\vspace{1cm} 
\justify }

\begin{document}

\newgeometry{top=25mm,bottom=25mm,left=25mm,right=25mm}
\maketitle
\thispagestyle{empty}
\Dedication

\begin{abstract}
\noindent
The organisers of major sports competitions use different policies with respect to constraints in the group draw. Our paper aims to rationalise these choices by analysing the trade-off between attractiveness (the number of games played by teams from the same geographic zone) and fairness (the departure of the draw mechanism from a uniform distribution). A parametric optimisation model is formulated and applied to the 2018 and 2022 FIFA World Cup draws. A flaw of the draw procedure is identified: the pre-assignment of the host to a group unnecessarily increases the distortions. All Pareto efficient sets of draw constraints are determined via simulations. The proposed framework can be used to find the optimal draw rules and justify the non-uniformity of the draw procedure for the stakeholders.
\end{abstract}


\AMS{62-08, 90-10, 90B90, 91B14}
% Computational methods for problems pertaining to statistics
% Mathematical modeling or simulation for problems pertaining to operations research and mathematical programming
% Case-oriented studies in operations research
% Social Choice

\JEL{C44, C63, Z20}
% Operations Research, Statistical Decision Theory
% Computational Techniques, Simulation Modeling 
% Sports Economics, General

\clearpage
\restoregeometry

\section{Introduction} \label{Sec1}

Several sports tournaments start with a group stage: the teams are allocated into groups to play a round-robin contest---where each team plays the same number of games against every other team---and the best teams in each group qualify for the next stage. The process of assigning the teams to the groups is called the \emph{group draw}.
The group draws of major sports competitions are among the most followed randomisation procedures in practice. For example, the video of the 2022 FIFA World Cup draw (\url{https://www.youtube.com/watch?v=PgV_z-TA0aw})---the draw that will serve as a case study for our optimisation model---has been watched 4.9 million times as of 30 September 2025.

According to recent papers in operational research (see Section~\ref{Sec2}) and the policies of governing bodies in sports, the most important requirements for a group draw are \emph{balancedness}, \emph{transparency}, and \emph{uniform distribution} \citep{Guyon2015a}.
Balance means that the groups must be roughly at the same competitive level; the probability of qualification is not allowed to depend substantially on the outcome of the draw.
Transparency is connected to the implementation of the draw; it typically involves celebrities drawing balls from a small number of urns and is broadcast live.
Uniform distribution is a fairness principle: every possible draw outcome should have the same probability of occurring.

In an unconstrained draw, the three criteria can be easily guaranteed.
Balancedness is usually achieved by seeding the teams into pots based on an exogenous ranking determined by historical performances, and creating groups that contain (at most) one team from each pot. This allows for a nice television show: the balls representing the teams in a pot are placed in an urn, and they are drawn sequentially to provide a sense of excitement and anticipation, without threatening transparency and uniform distribution.

However, the set of teams playing in the same group may be subject to further restrictions. In particular, \emph{geographic separation} is an important issue in the Olympic Games and World Championships because the stakeholders want to see unique matches between teams coming from different geographic regions. If geographical constraints need to be met, the three basic criteria may be violated since no mechanism is known that can satisfy all of them under arbitrary restrictions. The F\'ed\'eration Internationale de Football Association (FIFA), the international governing body of (association) football, has often sacrificed balancedness in the history of the FIFA World Cup \citep{Jones1990, RathgeberRathgeber2007, Guyon2015a}.
After \citet{Guyon2015a} identified several flaws in the 2014 FIFA World Cup draw, FIFA adopted the so-called \emph{Skip mechanism} for both the 2018 and 2022 FIFA World Cup draws \citep{Guyon2018d, Csato2023d}.

The Skip mechanism ensures that the teams can be drawn from urns associated with the seeding pots, while background computer calculations---that remain deterministic and could be verified by (sophisticated) viewers---guarantee geographic separation. Consequently, balancedness and transparency still hold, but uniform distribution might be violated \citep{RobertsRosenthal2024}. Non-uniformity is not merely of theoretical interest; the draw procedure can substantially distort the probability of qualification for some teams \citep{Csato2025c}.
Indeed, an indirect proof of the importance of non-uniformity is that, in contrast to the recent FIBA Basketball World Cups \citep{FIBA2019, FIBA2023} and FIFA World Cups \citep{FIFA2017c, FIFA2022a, FIFA2025b}, the 2025 IHF Men's World Handball Championship \citep{IHF2025} and the 2025 FIVB Men's Volleyball World Championship \citep{Volleyball2024} do not require geographic separation in their group draw.

These distinct rules highlight the trade-off between geographic separation and uniform distribution in the group draw of a sports tournament. However, the previous literature has not investigated this issue, which motivates the current study.
Our main contributions can be summarised as follows:
\begin{itemize}
\item
A parametric optimisation model is developed to find the best compromise between fairness (measured by the distortion compared to a uniform draw) and geographical separation (measured by the average number of intra-continental games) in the group draw of a sports tournament (Section~\ref{Sec3});
\item
The model is applied to the 2018 and 2022 FIFA World Cup draws, which involve extensive numerical simulations for 32 reasonable sets of draw restrictions (Section~\ref{Sec4});
\item
The official policy of pre-assigning the host to the first group is demonstrated to be suboptimal with respect to the bias of the draw (Section~\ref{Sec41});
\item
Depending on the composition of the draw pots and the measure of non-uniformity, at most five sets of constraints could be optimal according to the preferences of the organiser (Section~\ref{Sec43}).
\end{itemize}
Thus, we have a clear policy recommendation to improve the FIFA World Cup draw, as well as the draw of several major sports tournaments. 
The suggested approach can also be directly used to optimise the set of draw constraints before the draw takes place, and justify the unavoidable distortion of the transparent draw procedure for all stakeholders.
%To summarise, our results make it possible to retain the transparency of a group draw by uncovering the extent of the implied distortions.

\section{Related literature} \label{Sec2}

As we have seen in Section~\ref{Sec1}, balancedness, transparency, and uniform distribution can be easily guaranteed in any group draw in the absence of geographical separation. Hence, the interesting case is when geographic restrictions apply. Since transparency cannot be sacrificed for the sake of other goals in practice, FIFA has struggled with balancedness and uniform distribution in the history of the FIFA World Cup draw. The draws of the 1990 \citep{Jones1990}, 2006 \citep{RathgeberRathgeber2007}, and 2014 \citep{Guyon2015a} events were revealed to violate uniform distribution. The 2014 FIFA World Cup draw was quite imbalanced, too \citep{Guyon2015a}.
Therefore, \citet{Guyon2014a}---the full preliminary version of the published paper \citet{Guyon2015a}---has proposed three fairer and still transparent draw procedures, which inspired FIFA to adopt the Skip mechanism for the 2018 and 2022 FIFA World Cup draws \citep{Guyon2018d}.
\citet{Tijms2015} underlines the importance of transparency by discussing the credibility of accusations that the draw for the quarterfinals in the 2012/13 UEFA Champions League was manipulated.

\citet{LaprePalazzolo2023} empirically assess competitive imbalance between the FIFA World Cup groups from 1954 to 2022 to show its serious consequences with respect to the chance of reaching the quarterfinals. Consequently, reducing imbalance should be an important aim for the organiser.
\citet{CeaDuranGuajardoSureSiebertZamorano2020} develop a mixed integer linear programming model to improve the balance between the groups with retaining geographical separation.
\citet{LalienaLopez2019} propose two draw systems for allocating 32 teams into eight balanced groups subject to arbitrary constraints. The two mechanisms represent the trade-off between balancedness and transparency.
\citet{Csato2023d} demonstrates that the official seeding policy of the 2022 FIFA World Cup draw has increased imbalance between the groups; a better alternative would have been assigning the placeholders according to the highest-ranked potential winner rather than automatically seeding them in the weakest pot (see Section~\ref{Sec32}).
The methods of \citet{LalienaLopez2025} for balancing the strengths of groups containing three teams start from the ideas of \citet{LalienaLopez2019} but contain exogenous lower and upper bounds. The authors also examine non-uniformity but do not consider changing the geographic criteria.

Two published papers study the Skip mechanism.
\citet{RobertsRosenthal2024} explain how it can be simulated and compute the bias of the 2022 FIFA World Cup draw. Some procedures are presented to achieve uniform distribution by compromising transparency, as they contain some hidden randomness.
\citet{Csato2025c} investigates the effect of the draw order on the distortion of the 2018 FIFA World Cup draw. The pot sequence of the official method (Pot 1, Pot 2, Pot 3, Pot 4) turns out to be optimal, although even the small biases can substantially affect qualification probabilities.

However, none of these works discusses the possible gains from slacking the geographic constraints, even though these constraints have never been used in handball and volleyball. According to our knowledge, only one result exists on this issue in the literature.
The scope for improving fairness distortions in the UEFA Champions League Round of 16 draw \citep{KlossnerBecker2013} is found to be limited through a better draw procedure but a quantitatively large reduction can be achieved by allowing one match (instead of zero) between two teams from the same country \citep[Section~6.2]{BoczonWilson2023}.
Nonetheless, there are fundamental differences:
(1) the fairness distortion measure of \citet{BoczonWilson2023} does not quantify the level of non-uniformity (it can exceed zero even if the draw procedure has a uniform distribution over all valid assignments); 
(2) the setting of a bipartite graph with 16 nodes (UEFA Champions League Round of 16 draw) is much simpler than the eight groups of four teams in the 2022 FIFA World Cup draw; and
(3) the UEFA Champions League Round of 16 draw has used a distinct mechanism \citep{Csato2025f}.

The new incomplete round-robin format of UEFA club competitions also calls for a draw. Since the 2024/25 season, the UEFA Champions League starts with a league phase, where 36 teams play against eight opponents each and are ranked in a single table. \citet{DevriesereGoossensSpieksma2025} present an integer program to check whether the schedule can be completed if a match is drawn.
\citet{GuyonBenSalemBuchholtzerTanre2025} investigate four different draw procedures for the league phase. They are implemented using integer programming to express the draw constraints. The authors compare these mechanisms by Monte Carlo simulations to evaluate their fairness and compute the matchup probabilities---what we will do for the 2018 and 2022 FIFA World Cups under 32 reasonable sets of draw constraints.

The current paper is related to the numerous papers investigating various design issues in the FIFA World Cup.
\citet{MonksHusch2009} estimate the effect of home continent and seeding on the probability of qualification.
\citet{ScarfYusof2011} explore the impact of the seeding policy without considering geographical constraints.
%\citet{Csato2023c} shows that the probability of qualifying for the FIFA World Cup is not independent of the confederation to which a particular team is assigned; for example, a South American team could have tripled its chances by playing in Asia.
\citet{Guyon2020a}, \citet{ChaterArrondelGayantLaslier2021}, and \citet{Stronka2024} examine the connection between the schedule of group matches and their competitiveness.
\citet{Stronka2020} focus on the temptation to deliberately lose a match in the last round of the group stage in order to meet a preferred opponent in the following knockout stage. % and \citet{Csato2025b}
\citet{KrumerMoreno-Ternero2023} suggest fair allocations of qualifying berths among the continents. % and \citet{CsatoKissSzadoczki2025}
\citet{Guyon2020a}, \citet{GuajardoKrumer2024}, and \citet{CsatoGyimesi2026} propose alternative formats for the 2026 FIFA World Cup, the first event with 48 teams.
A recent survey of these studies can be found in \citet{DevriesereCsatoGoossens2025}.

Finally, we develop an optimisation problem that allows the organiser to maximise a weighted combination of attractiveness (in our case, more inter-continental group matches) and fairness (here, smaller distortion compared to a uniform draw). This is a standard approach in operational research applied to sports problems. \citet{RibeiroUrrutiadeWerra2023b} construct tournament schedules to minimise the number of breaks, the number of rounds where breaks occur, the difference between the number of breaks for the teams, and the carryover effect.
\citet{CsatoMolontayPinter2024} compare 12 valid schedules for the UEFA Champions League group stage with respect to the probability of weakly (where one team is indifferent) and strongly (where both teams are indifferent) stakeless matches.
\citet{VanBulckPaakkonenJacquetGoossens2026} propose a heuristic bilevel optimisation method to ensure a balanced emphasis on physical abilities and cognitive skills in rogaining, an orienteering running sport.
\citet{DaltropheShrotAronshtam2026} study the complexity of manipulating knockout tournaments for various types of manipulation, such as in favour of a single player achieving the highest score, or in favour of a given set of players achieving a higher score than the complementary set.

\section{Methodology} \label{Sec3}

This section is structured as follows.
Section~\ref{Sec31} discusses the Skip mechanism and the implementation of a uniform draw.
Section~\ref{Sec32} outlines the 2022 FIFA World Cup draw, which motivates the optimisation problems developed in Section~\ref{Sec33}.
Finally, Section~\ref{Sec34} presents our simulation procedure.

\subsection{Draw procedures} \label{Sec31}

Draw constraints such as geographic separation make an assignment of the teams to the groups \emph{invalid} if it violates at least one constraint. The statistical challenge resides in simulating (closely) from the uniform distribution with a sequential and entertaining method. Governing bodies in sports use two procedures for this purpose \citep{Csato2025f}.
Among them, the Skip mechanism is more popular as it is employed in the FIBA Basketball World Cup, the FIFA World Cup, and the UEFA Nations League, among other tournaments.
Its underlying idea is straightforward: if the violation of a draw constraint occurs or is anticipated to occur, the computer used to assist with the draw indicates the next available group in alphabetical order for the drawn team \citep{UEFA2024a}.

However, while direct conflicts are straightforward to detect, the set of groups available to a team might depend not only on its attributes and those of the previously drawn teams, but also on the attributes of the teams still to be drawn. In other words, assigning a team to a certain group may violate no constraint immediately, but would make it impossible to assign the remaining teams to the free slots. Hence, the computer calculations should anticipate every possible scenario and prevent any deadlock situation, that is, at least one valid solution needs to remain for the teams still to be drawn as illustrated by the following example.

\begin{example} \label{Examp1}
Assume that six teams 1--6 are divided into three pots of two teams to create two groups, A and B, each containing one team from each pot: Pot 1 consists of teams 1 and 2, Pot 2 consists of teams 3 and 4, and Pot 3 consists of teams 5 and 6.
There is a further draw constraint: at most two teams from the set of three teams $\{ 2, 4, 6 \}$ can play in the same group.

If the Skip mechanism starts with Pot 1 and finishes with Pot 3, the first conflict may arise only when the draw reaches Pot 3. If Group A contains teams 1 and 3, while Group B contains teams 2 and 4, two possibilities can occur:
\begin{itemize}
\item
Team 5 is drawn first from Pot 3 and is assigned to Group B by the Skip mechanism, even though Group A has a free slot for a team from Pot 3 (otherwise, a deadlock will occur as assigning team 6 to Group B violates the draw constraint); or
\item
Team 6 is drawn first from Pot 3 and is assigned to Group A.
\end{itemize}
\end{example}

Naturally, the identification of a deadlock is not always straightforward with several pots and teams, and simulating the Skip mechanism with a computer program is surprisingly challenging \citep{RobertsRosenthal2024}. According to our knowledge, no governing body publishes the exact algorithm used for the draw, but appropriate backtracking algorithms have been proposed in the academic literature \citep{Csato2025c, Guyon2014a, RobertsRosenthal2024}.

%Example~\ref{Examp2} also reveals that the Skip mechanism has $k!$ variants if the number of pots is $k$.

Obviously, the Skip mechanism always leads to a feasible solution if one exists. However, it is not guaranteed that each valid assignment has the same probability of occurring, the Skip mechanism is known to be non-uniform \citep{Csato2025c, RobertsRosenthal2024} as in Example~\ref{Examp1}.
Here, there are four possible outcomes if the group labels are ignored, but one of them (when teams 2, 4, 6 are assigned to the same group) is invalid, resulting in three feasible assignments. In particular, team 1 can be in the same group with
(a) teams 3 and 6;
(b) teams 4 and 5; or
(c) teams 4 and 6.
Thus, teams 1 and 4 play in the same group with a probability of 2/3 in a uniform draw.
The Skip mechanism does not consider the draw constraint until Pot 3 is reached; therefore, teams 1 and 4 play in the same group with a probability of 1/2.

A uniform draw can be implemented, for instance, by a rejection sampler \citep[Section~2.1]{RobertsRosenthal2024}. This works as follows. First, an unconstrained draw is generated: the teams are drawn sequentially and selected uniformly at random from the pots. Once the draw outcome is determined, it is checked whether all draw constraints are satisfied. If yes, the output is the chosen draw. If not, a new unconstrained draw is generated.
Obviously, the rejection sampler is not used in the draw of sports tournaments because the number of repetitions required to generate a valid assignment is highly volatile.
%Luckily, generating one unconstrained draw can be sufficient. But in the case of the 2022 FIFA World Cup draw, the probability of finding a feasible solution is only about 1/ (Section~\ref{Sec4}), therefore, sometimes thousands of runs are needed to get the desired assignment, which satisfies all constraints.

\subsection{The rules of the 2018 and 2022 FIFA World Cup draws} \label{Sec32}

\citet{Guyon2015a} identified a number of flaws in the 2014 FIFA World Cup draw, and managed to convince FIFA to reform the draw of the FIFA World Cup \citep{Guyon2018d}. 
Therefore, first in the history, the draw pots were constructed by considering the historical performances of the national teams, while geographic separation was achieved by draw constraints in both the 2018 and 2022 tournaments.

In particular, the 32 national teams were divided into four pots of eight teams each. Pot 1 contained the host, automatically assigned to Group A, as well as the seven highest-ranked teams in the FIFA World Rankings. Pot 2 consisted of the next eight highest-ranked teams, and Pot 3 consisted of the next eight. Finally, Pot 4 contained the eight lowest-ranked nations in 2018, and the five lowest-ranked nations and three placeholders---the two winners of the inter-confederation play-offs and one winner of a UEFA play-off still to be played---in 2022.
One team was assigned from each pot to a group to ensure roughly the same competitive level in all groups.

FIFA has six confederations: AFC (Asian Football Confederation), CAF (Confederation of African Football), CONCACAF (Confederation of North, Central American and Caribbean Association Football), CONMEBOL (Confederación Sudamericana de Fútbol), OFC (Oceania Football Confederation), and UEFA (Union of European Football Associations).
In the qualification for the FIFA World Cup, each confederation received a predetermined number of slots, possibly including shared slots allocated by inter-confederation play-offs. In order to guarantee geographic separation and maximise the number of inter-continental games in the group stage, the following restrictions were used \citep{FIFA2017c, FIFA2022a}:
\begin{itemize}
\item
\emph{Constraint A}: No group can contain more than one AFC team;
\item
\emph{Constraint B}: No group can contain more than one CAF team;
\item
\emph{Constraint C}: No group can contain more than one CONCACAF team;
\item
\emph{Constraint D}: No group can contain more than one CONMEBOL team;
\item
\emph{Constraint E}: No group can contain fewer than one or more than two UEFA teams.
\end{itemize}

\begin{table}[t!]
  \centering
  \caption{Matches played in the 2022 FIFA World Cup qualification after the draw}
  \label{Table1}
\begin{threeparttable}
  \rowcolors{3}{}{gray!20}
    \begin{tabularx}{\textwidth}{l lcLc} \toprule
    Match & Team 1 & Elo rating 1 & Team 2 & Elo rating 2 \\ \bottomrule
    AFC PO & Australia & 1677  & United Arab Emirates & 1515 \\
    IPO1 & Peru (CONMEBOL) & 1856  & Winner of AFC PO & 1677/1515 \\ \hline
    IPO2 & Costa Rica (CONC) & 1743  & New Zealand (OFC) & 1557 \\ \toprule
    \end{tabularx}
\begin{tablenotes} \footnotesize
\item
All play-offs are single-leg games played on a neutral ground in Qatar.
\item
CONC stands for CONCACAF.
\item
Elo ratings show World Football Elo Ratings on 31 March 2022 (the day of the draw), see \url{https://www.international-football.net/elo-ratings-table?year=2022&month=03&day=31}.
\end{tablenotes}
\end{threeparttable}
\end{table}

In the 2022 FIFA World Cup draw, the placeholders for the inter-confederation play-offs were considered for all confederations involved. The qualification games played after the draw are shown in Table~\ref{Table1}.

\begin{table}[t!]
  \centering
  \caption{The geographical composition of the pots in the recent FIFA World Cups}
  \label{Table2}

\begin{subtable}{\textwidth}
  \centering
  \caption{2018 FIFA World Cup draw}
  \label{Table2a}
  \rowcolors{3}{}{gray!20}
    \begin{tabularx}{0.8\textwidth}{l CCCCC} \toprule
    Confederation & Pot 1 & Pot 2 & Pot 3 & Pot 4 & Sum \\ \bottomrule
    AFC   & 1     & 0     & 1     & 4     & 5 \\
    CAF   & 0     & 0     & 3     & 2     & 5 \\
    CONCACAF & 0     & 1     & 1     & 1     & 3 \\
    CONMEBOL & 2     & 3     & 0     & 0     & 5 \\
    UEFA  & 6     & 4     & 3     & 1     & 14 \\ \toprule
    \end{tabularx}
\end{subtable}

\vspace{0.25cm}
\begin{subtable}{\textwidth}
  \centering
  \caption{2022 FIFA World Cup draw}
  \label{Table2b}
  
\begin{threeparttable}
  \rowcolors{3}{}{gray!20}
    \begin{tabularx}{0.8\textwidth}{l CCCCC} \toprule
    Confederation & Pot 1 & Pot 2 & Pot 3 & Pot 4 & Sum \\ \bottomrule
    AFC   & 1     & 0     & 3     & 2     & 6 \\
    CAF   & 0     & 0     & 3     & 2     & 5 \\
    CONCACAF & 0     & 2     & 0     & 2     & 4 \\
    CONMEBOL & 2     & 1     & 0     & 2     & 5 \\
    OFC   & 0     & 0     & 0     & 1     & 1 \\
    UEFA  & 5     & 5     & 2     & 1     & 13 \\ \bottomrule
    \end{tabularx}
\begin{tablenotes} \footnotesize
\item
Placeholders for the inter-confederation play-offs are taken for both confederations into account, see Table~\ref{Table1}. 
\end{tablenotes}
\end{threeparttable}
\end{subtable}

\end{table}

Table~\ref{Table2} shows the number of teams in the four pots from the six confederations.
In 2018 (Table~\ref{Table2a}), no OFC team qualified for the FIFA World Cup. In 2022 (Table~\ref{Table2b}), at most one OFC team could have played in the group stage; hence, there was no need to avoid a match between two OFC teams by an additional restriction.

It is worth noting that the 2026 FIFA World Cup draw also used Constraints A--E and the Skip mechanism, with pre-assigning the three hosts. However, each pot contained 12 teams instead of eight, and additional restrictions were imposed for the four strongest teams in Pot 1 to ensure competitive balance \citep{FIFA2025b}.

\subsection{The optimisation model} \label{Sec33}

\begin{table}[t!]
  \centering
  \caption{Sets of possible draw constraints in the recent FIFA World Cup draws}
  \label{Table3}
  \rowcolors{3}{}{gray!20}
\centerline{
    \begin{tabularx}{1.15\textwidth}{l CCCCC} \toprule \hiderowcolors
    Confederation & AFC & CAF  & CONCACAF & CONMEBOL & UEFA \\
    Teams in a group & At most 1 & At most 1 & At most 1 & At most 1 & 1 or 2 \\
    Scenario & Constraint A & Constraint B & Constraint C & Constraint D & Constraint E \\ \bottomrule \showrowcolors
    0     & \textcolor{BrickRed}{\ding{55}} & \textcolor{BrickRed}{\ding{55}} & \textcolor{BrickRed}{\ding{55}} & \textcolor{BrickRed}{\ding{55}} & \textcolor{BrickRed}{\ding{55}} \\
    1     & \textcolor{BrickRed}{\ding{55}} & \textcolor{BrickRed}{\ding{55}} & \textcolor{BrickRed}{\ding{55}} & \textcolor{BrickRed}{\ding{55}} & \textcolor{PineGreen}{\ding{52}} \\
    2     & \textcolor{BrickRed}{\ding{55}} & \textcolor{BrickRed}{\ding{55}} & \textcolor{BrickRed}{\ding{55}} & \textcolor{PineGreen}{\ding{52}} & \textcolor{BrickRed}{\ding{55}} \\
    3     & \textcolor{BrickRed}{\ding{55}} & \textcolor{BrickRed}{\ding{55}} & \textcolor{BrickRed}{\ding{55}} & \textcolor{PineGreen}{\ding{52}} & \textcolor{PineGreen}{\ding{52}} \\
    4     & \textcolor{BrickRed}{\ding{55}} & \textcolor{BrickRed}{\ding{55}} & \textcolor{PineGreen}{\ding{52}} & \textcolor{BrickRed}{\ding{55}} & \textcolor{BrickRed}{\ding{55}} \\
    5     & \textcolor{BrickRed}{\ding{55}} & \textcolor{BrickRed}{\ding{55}} & \textcolor{PineGreen}{\ding{52}} & \textcolor{BrickRed}{\ding{55}} & \textcolor{PineGreen}{\ding{52}} \\
    6     & \textcolor{BrickRed}{\ding{55}} & \textcolor{BrickRed}{\ding{55}} & \textcolor{PineGreen}{\ding{52}} & \textcolor{PineGreen}{\ding{52}} & \textcolor{BrickRed}{\ding{55}} \\
    7     & \textcolor{BrickRed}{\ding{55}} & \textcolor{BrickRed}{\ding{55}} & \textcolor{PineGreen}{\ding{52}} & \textcolor{PineGreen}{\ding{52}} & \textcolor{PineGreen}{\ding{52}} \\
    8     & \textcolor{BrickRed}{\ding{55}} & \textcolor{PineGreen}{\ding{52}} & \textcolor{BrickRed}{\ding{55}} & \textcolor{BrickRed}{\ding{55}} & \textcolor{BrickRed}{\ding{55}} \\
    9     & \textcolor{BrickRed}{\ding{55}} & \textcolor{PineGreen}{\ding{52}} & \textcolor{BrickRed}{\ding{55}} & \textcolor{BrickRed}{\ding{55}} & \textcolor{PineGreen}{\ding{52}} \\
    10    & \textcolor{BrickRed}{\ding{55}} & \textcolor{PineGreen}{\ding{52}} & \textcolor{BrickRed}{\ding{55}} & \textcolor{PineGreen}{\ding{52}} & \textcolor{BrickRed}{\ding{55}} \\
    11    & \textcolor{BrickRed}{\ding{55}} & \textcolor{PineGreen}{\ding{52}} & \textcolor{BrickRed}{\ding{55}} & \textcolor{PineGreen}{\ding{52}} & \textcolor{PineGreen}{\ding{52}} \\
    12    & \textcolor{BrickRed}{\ding{55}} & \textcolor{PineGreen}{\ding{52}} & \textcolor{PineGreen}{\ding{52}} & \textcolor{BrickRed}{\ding{55}} & \textcolor{BrickRed}{\ding{55}} \\
    13    & \textcolor{BrickRed}{\ding{55}} & \textcolor{PineGreen}{\ding{52}} & \textcolor{PineGreen}{\ding{52}} & \textcolor{BrickRed}{\ding{55}} & \textcolor{PineGreen}{\ding{52}} \\
    14    & \textcolor{BrickRed}{\ding{55}} & \textcolor{PineGreen}{\ding{52}} & \textcolor{PineGreen}{\ding{52}} & \textcolor{PineGreen}{\ding{52}} & \textcolor{BrickRed}{\ding{55}} \\
    15    & \textcolor{BrickRed}{\ding{55}} & \textcolor{PineGreen}{\ding{52}} & \textcolor{PineGreen}{\ding{52}} & \textcolor{PineGreen}{\ding{52}} & \textcolor{PineGreen}{\ding{52}} \\
    16    & \textcolor{PineGreen}{\ding{52}} & \textcolor{BrickRed}{\ding{55}} & \textcolor{BrickRed}{\ding{55}} & \textcolor{BrickRed}{\ding{55}} & \textcolor{BrickRed}{\ding{55}} \\
    17    & \textcolor{PineGreen}{\ding{52}} & \textcolor{BrickRed}{\ding{55}} & \textcolor{BrickRed}{\ding{55}} & \textcolor{BrickRed}{\ding{55}} & \textcolor{PineGreen}{\ding{52}} \\
    18    & \textcolor{PineGreen}{\ding{52}} & \textcolor{BrickRed}{\ding{55}} & \textcolor{BrickRed}{\ding{55}} & \textcolor{PineGreen}{\ding{52}} & \textcolor{BrickRed}{\ding{55}} \\
    19    & \textcolor{PineGreen}{\ding{52}} & \textcolor{BrickRed}{\ding{55}} & \textcolor{BrickRed}{\ding{55}} & \textcolor{PineGreen}{\ding{52}} & \textcolor{PineGreen}{\ding{52}} \\
    20    & \textcolor{PineGreen}{\ding{52}} & \textcolor{BrickRed}{\ding{55}} & \textcolor{PineGreen}{\ding{52}} & \textcolor{BrickRed}{\ding{55}} & \textcolor{BrickRed}{\ding{55}} \\
    21    & \textcolor{PineGreen}{\ding{52}} & \textcolor{BrickRed}{\ding{55}} & \textcolor{PineGreen}{\ding{52}} & \textcolor{BrickRed}{\ding{55}} & \textcolor{PineGreen}{\ding{52}} \\
    22    & \textcolor{PineGreen}{\ding{52}} & \textcolor{BrickRed}{\ding{55}} & \textcolor{PineGreen}{\ding{52}} & \textcolor{PineGreen}{\ding{52}} & \textcolor{BrickRed}{\ding{55}} \\
    23    & \textcolor{PineGreen}{\ding{52}} & \textcolor{BrickRed}{\ding{55}} & \textcolor{PineGreen}{\ding{52}} & \textcolor{PineGreen}{\ding{52}} & \textcolor{PineGreen}{\ding{52}} \\
    24    & \textcolor{PineGreen}{\ding{52}} & \textcolor{PineGreen}{\ding{52}} & \textcolor{BrickRed}{\ding{55}} & \textcolor{BrickRed}{\ding{55}} & \textcolor{BrickRed}{\ding{55}} \\
    25    & \textcolor{PineGreen}{\ding{52}} & \textcolor{PineGreen}{\ding{52}} & \textcolor{BrickRed}{\ding{55}} & \textcolor{BrickRed}{\ding{55}} & \textcolor{PineGreen}{\ding{52}} \\
    26    & \textcolor{PineGreen}{\ding{52}} & \textcolor{PineGreen}{\ding{52}} & \textcolor{BrickRed}{\ding{55}} & \textcolor{PineGreen}{\ding{52}} & \textcolor{BrickRed}{\ding{55}} \\
    27    & \textcolor{PineGreen}{\ding{52}} & \textcolor{PineGreen}{\ding{52}} & \textcolor{BrickRed}{\ding{55}} & \textcolor{PineGreen}{\ding{52}} & \textcolor{PineGreen}{\ding{52}} \\
    28    & \textcolor{PineGreen}{\ding{52}} & \textcolor{PineGreen}{\ding{52}} & \textcolor{PineGreen}{\ding{52}} & \textcolor{BrickRed}{\ding{55}} & \textcolor{BrickRed}{\ding{55}} \\
    29    & \textcolor{PineGreen}{\ding{52}} & \textcolor{PineGreen}{\ding{52}} & \textcolor{PineGreen}{\ding{52}} & \textcolor{BrickRed}{\ding{55}} & \textcolor{PineGreen}{\ding{52}} \\
    30    & \textcolor{PineGreen}{\ding{52}} & \textcolor{PineGreen}{\ding{52}} & \textcolor{PineGreen}{\ding{52}} & \textcolor{PineGreen}{\ding{52}} & \textcolor{BrickRed}{\ding{55}} \\
    31    & \textcolor{PineGreen}{\ding{52}} & \textcolor{PineGreen}{\ding{52}} & \textcolor{PineGreen}{\ding{52}} & \textcolor{PineGreen}{\ding{52}} & \textcolor{PineGreen}{\ding{52}} \\ \bottomrule
    \end{tabularx}
}
\end{table}

Both the 2018 and 2022 FIFA World Cup draws used the Skip mechanism, starting with Pot 1 and ending with Pot 4.
According to Section~\ref{Sec31}, this procedure is non-uniform, which might threaten the fairness of the draw. In the absence of a uniform and transparent draw procedure \citep{RobertsRosenthal2024}, a reasonable step towards uniformity is offered by ignoring some of Constraints A--E.
In particular, Constraints A to E can be separately imposed or disregarded, implying 32 different sets of draw constraints $\mathcal{C}$.
They are listed in Table~\ref{Table3} together with their notation by integers from 0 to 31.

Ignoring a constraint is expected to increase the number of inter-confederation matches played in the group stage, which reduces attractiveness since national teams from the same confederation play several matches outside the FIFA World Cup in its qualifiers and in continental championships.

Therefore, we assume that FIFA has two reasonable goals when designing the rules of the draw:
\begin{itemize}
\item
Decrease the distortion of the draw procedure;
\item
Maximise the number of inter-confederation group matches.
\end{itemize}
Since the number of group matches is fixed by the tournament format, the second objective is equivalent to minimising the number of intra-confederation games, which will be considered in the following.

The non-uniformity of the draw under a set of constraints $\mathcal{C}$ can be quantified in several ways; we will use two measures.
The first is the mean absolute bias of all non-zero matchup probabilities, measured in percentage points:
\begin{equation}
\Delta (\mathcal{C}) = 100 \cdot \frac{\sum_{i,j} \left| p_{ij}^U(\mathcal{C}) - p_{ij}^S(\mathcal{C}) \right|}{\sum_{i,j} \# \left\{ p_{ij}^U(\mathcal{C}) > 0 \right\}},
\end{equation}
where $p_{ij}^U(\mathcal{C})$ and $p_{ij}^S(\mathcal{C})$ are the probabilities that teams $i$ and $j$ play in the same group under a set of constraints $\mathcal{C}$ with a uniform draw ($U$) and the Skip mechanism ($S$), respectively, while $\sum_{i,j} \# \left\{ p_{ij}^U(\mathcal{C}) > 0 \right\}$ is the number of team pairs with a non-zero probability of being drawn into the same group.

The second is the maximum absolute bias of all matchup probabilities, measured in percentage points:
\begin{equation}
\Omega (\mathcal{C}) = 100 \cdot \max_{i,j} \left| p_{ij}^U(\mathcal{C}) - p_{ij}^S(\mathcal{C}) \right|.
\end{equation}

These two concepts are analogous to the \emph{average case} ($\Delta$) and the \emph{worst case} ($\Omega$) for the resource usage of an algorithm in computer science.

On the other hand, only one metric is considered for the second goal, the average number of intra-confederation group games $\Psi (\mathcal{C})$.
Note that the maximal number of intra-confederation games would not be informative without knowing the probability of such an unfavourable outcome, in contrast to the maximum absolute bias, which is certainly relevant to at least one pair of teams.

Our measures for the two goals of the organiser imply two optimisation problems:
\begin{equation} \label{Opt1}
\min_\mathcal{C} \alpha \cdot 10 \Delta (\mathcal{C}) + (1 - \alpha) \cdot \Psi (\mathcal{C}),
\end{equation}
\begin{equation} \label{Opt2}
\min_\mathcal{C} \alpha \cdot \Omega (\mathcal{C}) + (1 - \alpha) \cdot \Psi (\mathcal{C}),
\end{equation}
where parameter $\alpha$ reflects the priority of the decision-maker with respect to the objectives. In other words, the aim is to minimise the weighted average of non-uniformity and intra-confederation games. Multiplier 10 is added to $\Delta (\mathcal{C})$ in~\eqref{Opt1} based on our two case studies, as the values of $10 \Delta (\mathcal{C})$ and $\Psi (\mathcal{C})$ are roughly of the same magnitude.

We will determine the set of \emph{Pareto efficient draw constraints} that can be the optimal solutions to optimisation problems~\eqref{Opt1} and \eqref{Opt2} for a given value of parameter $\alpha$. Naturally, they might be different for the 2018 and 2022 FIFA World Cups, but the multiplier 10 in~\eqref{Opt1} does not affect them. The organiser could explain these choices to the public if the decision criteria are the non-uniformity of the draw and the number of intra-confederation group games.
It is straightforward to obtain the optimal solution for any exogenous value of $\alpha$ because 32 sets of draw constraints are considered in each case; thus, the objective function should be evaluated only in these 32 elements of set $\mathcal{C}$.

\subsection{Details of the computation} \label{Sec34}

Drawing the teams in random order, the number of possible outcomes in both the 2018 and 2022 FIFA World Cup draws is $\left( 8! \right)^3 \approx 6.6 \cdot 10^{13}$ with a uniform draw and $7! \cdot \left( 8! \right)^3 \approx 3.3 \cdot 10^{17}$ with the Skip mechanism. Therefore, they are simulated one million times, similar to \citet{Csato2025c} and \citet{RobertsRosenthal2024} (10 million is used for the results in Section~\ref{Sec41}). In the case of uniform draw, which is implemented by a rejection sampler \citep[Section~2.1]{RobertsRosenthal2024}, we generate one million valid assignments under all Constraints A--E (scenario 31 in Table~\ref{Table3}). However, this leads to a random sample with (much) more feasible solutions in scenarios 0--30 since the chosen unconstrained draw is not rejected if only the constraints disregarded in the given scenario are violated. For instance, we consider about 2.25 (4.7) million uniform draws for scenario 25 (Constraints A, B, E) in the 2018 (2022) FIFA World Cup draw.

For the Skip mechanism, the host is \emph{not} automatically placed in Group A since this may unnecessarily distort the matchup probabilities due to the sequential nature of the Skip mechanism, as will be demonstrated in Section~\ref{Sec41}. Note that, if the host is automatically placed in Group A, it will play against any team from Pot 2 with an equal probability of (exactly) 12.5\%. However, these probabilities are certainly different under a uniform draw because both Pots 1 and 2 contain CONMEBOL teams that cannot play against each other (Table~\ref{Table2}).
Since the groups could be labelled even \emph{after} the draw is finished, the group to which the host is assigned can be called Group A without any pre-assignment.

The probabilities $p_{ij}^U(\mathcal{C})$ and $p_{ij}^S(\mathcal{C})$ come directly from these simulated draws.
In the calculation of $\Psi (\mathcal{C})$ for the 2022 FIFA World Cup draw, the uncertainty of the winners of inter-confederation play-offs IPO1 and IPO2 should also be taken into account. They are determined according to the Elo ratings reported in Table~\ref{Table1}, following the win expectancy $W_E$ of World Football Elo Ratings (\url{https://eloratings.net/about}):
\[
W_{ij}^E = \frac{1}{1 + 10^{- \left( R_i - R_j \right)/400}}.
\]
Consequently,
\begin{itemize}
\item
the probability that Australia wins the AFC PO is 71.76\%;
\item
the probability that Peru wins the IPO1 is 77.65\%;
\item
the probability that Costa Rica wins the IPO2 is 74.47\%;
\end{itemize}
For example, if the particular draw assigns an AFC team and the winner of IPO1 to the same group, the number of group matches between two AFC teams is increased only by $0.2235 = 1 - 0.7765$ because Peru is not an AFC team.

\section{Results} \label{Sec4}

Our findings on the FIFA World Cup draw are presented in three parts, each of them having a key message for decision-makers.
Section~\ref{Sec41} shows the effect of pre-assigning the host to Group A and offers a better method without any drawback.
Section~\ref{Sec42} focuses on the restrictiveness and the non-uniformity of the draw under the 32 possible sets of constraints.
Last but not least, Section~\ref{Sec43} identifies the Pareto frontier with respect to the trade-off between the distortion of the Skip mechanism and the number of group matches played by teams from the same continent.

\subsection{The suboptimality of pre-assigning the host} \label{Sec41}

\begin{figure}[t!]
\centering

\begin{subfigure}{\textwidth}
\centering
\caption{2018 FIFA World Cup draw}
\label{Fig1a}

\begin{tikzpicture}
%\selectcolormodel{gray}
\begin{axis}[
name = axis1,
width = 0.8\textwidth, 
height = 0.64\textwidth,
xmajorgrids,
ymajorgrids,
%xbar stacked,
%bar width = 10pt,
scaled x ticks = false,
xlabel = {Mean absolute distortion (\%)},
xlabel style = {align=center, font=\small},
xticklabel style = {/pgf/number format/fixed,/pgf/number format/precision=5},
xmin = 0,
%extra x ticks = 0,
%extra x tick labels = ,
%extra x tick style = {grid = major, major grid style = {black,very thick}},
symbolic y coords = {Russia,Germany,Brazil,Portugal,Argentina,Belgium,Poland,France,Spain,Peru,Switzerland,England,Colombia,Mexico,Uruguay,Croatia,Denmark,Iceland,Costa Rica,Sweden,Tunisia,Egypt,Senegal,Iran,Serbia,Nigeria,Australia,Japan,Morocco,Panama,South Korea,Saudi Arabia},
ytick = data,
y dir = reverse,
enlarge y limits = 0.02,
]
% Our proposal: Groups labelled after the draw
\addplot [red, mark=pentagon, only marks, mark size=2pt, mark options={solid,semithick}] coordinates{
(0.446638333333333,Russia)
(0.447778333333333,Germany)
(1.52134761904762,Brazil)
(0.442146666666667,Portugal)
(1.51672380952381,Argentina)
(0.439164166666667,Belgium)
(0.4471925,Poland)
(0.447545833333334,France)
(0.681038333333333,Spain)
(0.467886363636364,Peru)
(0.679348333333333,Switzerland)
(0.680529166666667,England)
(0.466812727272727,Colombia)
(2.89802454545455,Mexico)
(0.46904909090909,Uruguay)
(0.67668,Croatia)
(0.437815,Denmark)
(0.4374,Iceland)
(0.296624545454545,Costa Rica)
(0.432098333333334,Sweden)
(0.253389090909091,Tunisia)
(0.256443636363636,Egypt)
(0.2502,Senegal)
(0.593973,Iran)
(2.0920475,Serbia)
(0.284831428571429,Nigeria)
(0.344303478260869,Australia)
(0.341070434782609,Japan)
(0.278688571428571,Morocco)
(0.564488181818182,Panama)
(0.347185217391305,South Korea)
(0.346874782608696,Saudi Arabia)
};
% Official rule: host is pre-assigned
\addplot [blue, mark=asterisk, only marks, mark size=2pt, mark options={solid,semithick}] coordinates{
(2.92669,Russia)
(0.5438125,Germany)
(1.70302857142857,Brazil)
(0.5465,Portugal)
(1.71211904761905,Argentina)
(0.547515,Belgium)
(0.545241666666667,Poland)
(0.545893333333334,France)
(0.610325833333333,Spain)
(0.7769,Peru)
(0.610603333333334,Switzerland)
(0.6160575,England)
(0.771305454545454,Colombia)
(2.84581363636364,Mexico)
(0.778398181818182,Uruguay)
(0.615999166666667,Croatia)
(0.673001666666666,Denmark)
(0.677005,Iceland)
(0.36932,Costa Rica)
(0.672179166666666,Sweden)
(0.475612727272727,Tunisia)
(0.474937272727273,Egypt)
(0.482447272727273,Senegal)
(0.774622,Iran)
(2.14012166666667,Serbia)
(0.40025619047619,Nigeria)
(0.381928695652174,Australia)
(0.384102608695652,Japan)
(0.399451428571429,Morocco)
(0.575505454545455,Panama)
(0.385164347826087,South Korea)
(0.382678260869565,Saudi Arabia)
};
\end{axis}
\end{tikzpicture}
\end{subfigure}

\vspace{0.25cm}
\begin{subfigure}{\textwidth}
\centering
\caption{2022 FIFA World Cup draw}
\label{Fig1b}
	
\begin{tikzpicture}
%\selectcolormodel{gray}
\begin{axis}[
name = axis1,
width = 0.8\textwidth, 
height = 0.64\textwidth,
xmajorgrids,
ymajorgrids,
%xbar stacked,
%bar width = 10pt,
scaled x ticks = false,
xlabel = {Mean absolute distortion (\%)},
xlabel style = {align=center, font=\small},
xticklabel style = {/pgf/number format/fixed,/pgf/number format/precision=5},
xmin = 0,
%extra x ticks = 0,
%extra x tick labels = ,
%extra x tick style = {grid = major, major grid style = {black,very thick}},
symbolic y coords = {Qatar,Brazil,Belgium,France,Argentina,England,Spain,Portugal,Mexico,Netherlands,Denmark,Germany,Uruguay,Switzerland,United States,Croatia,Senegal,Iran,Japan,Morocco,Serbia,Poland,South Korea,Tunisia,Cameroon,Canada,Ecuador,Saudi Arabia,Ghana,IPO1,IPO2,UEFA PO},
ytick = data,
y dir = reverse,
enlarge y limits = 0.02,
legend style = {font=\small,at={(0,-0.15)},anchor=north west,legend columns=1},
legend entries = {The group with the host is called Group A (our proposal)$\qquad$, The host is pre-assigned to Group A (official draw procedure)$\,$},
]
% Our proposal: Groups labelled after the draw
\addplot [red, mark=pentagon, only marks, mark size=2pt, mark options={solid,semithick}] coordinates{
(2.314,Qatar)
(2.49810571428571,Brazil)
(1.1176875,Belgium)
(1.11706333333333,France)
(2.4972,Argentina)
(1.11724666666667,England)
(1.118885,Spain)
(1.11013166666667,Portugal)
(2.41439272727273,Mexico)
(0.90219,Netherlands)
(0.909475833333333,Denmark)
(0.906444166666667,Germany)
(1.603959,Uruguay)
(0.908005,Switzerland)
(2.41671272727273,United States)
(0.906438333333333,Croatia)
(0.853497272727272,Senegal)
(1.01082857142857,Iran)
(1.01416285714286,Japan)
(0.855031818181818,Morocco)
(0.947915,Serbia)
(0.94799,Poland)
(1.02455047619048,South Korea)
(0.849442727272728,Tunisia)
(0.506946666666666,Cameroon)
(0.575732727272727,Canada)
(1.37674,Ecuador)
(1.556023,Saudi Arabia)
(0.501605714285714,Ghana)
(1.11745411764706,IPO1)
(0.570482727272727,IPO2)
(3.09464833333333,UEFA PO)
};
% Official rule: host is pre-assigned
\addplot [blue, mark=asterisk, only marks, mark size=2pt, mark options={solid,semithick}] coordinates{
(4.15904736842105,Qatar)
(2.45814380952381,Brazil)
(1.19548,Belgium)
(1.18858416666667,France)
(2.45700380952381,Argentina)
(1.19300833333333,England)
(1.18879,Spain)
(1.19056666666667,Portugal)
(2.45886636363636,Mexico)
(0.906918333333333,Netherlands)
(0.903209166666667,Denmark)
(0.9012925,Germany)
(1.434829,Uruguay)
(0.9037925,Switzerland)
(2.46208909090909,United States)
(0.902066666666667,Croatia)
(1.21285909090909,Senegal)
(1.12083619047619,Iran)
(1.12150095238095,Japan)
(1.20825909090909,Morocco)
(1.5214,Serbia)
(1.51341416666667,Poland)
(1.13069809523809,South Korea)
(1.20155909090909,Tunisia)
(0.601135238095238,Cameroon)
(0.533460909090909,Canada)
(1.91892380952381,Ecuador)
(1.873359,Saudi Arabia)
(0.597937142857143,Ghana)
(1.00162117647059,IPO1)
(0.533529090909091,IPO2)
(3.35127416666667,UEFA PO)
};
\end{axis}
\end{tikzpicture}
\end{subfigure}

\caption{The average absolute bias of the Skip mechanism for all countries}
\label{Fig1}

\end{figure}

%\end{document}

According to Section~\ref{Sec34}, the host can be guaranteed to play in Group A without pre-assigning it to the first group in alphabetical order as done by FIFA.
Figure~\ref{Fig1} shows the average absolute distortion of the Skip mechanism at the level of national teams for these two draw policies, as well as for both FIFA World Cups. In particular, the absolute biases are summed for all opponents of a given national team, and they are averaged. In both cases, the draw has the highest distortion for the host, Russia (2018) and Qatar (2022), although Mexico is also treated quite unfairly by the 2018 FIFA World Cup draw. As can be seen, our proposal reduces the distortion for the host by 85\% and 44\%, respectively, while it essentially does not worsen the situation of any other national team.
The distortion is also reduced for the majority of team pairs that can play in the same group: for 196 pairs out of the total 365 (53.7\%) in 2018 and for 233 pairs out of the total 355 (62.8\%) in 2022. Analogously, the mean absolute bias of the draw increases by 30\% (14\%) in the 2018 (2022) FIFA World Cup if the host is pre-assigned to Group A. 

These results are reinforced by Tables~\ref{Table_A1} and \ref{Table_A2} in the Appendix, which compare the distortions for all team pairs under the two versions of the Skip mechanism. The advantage of the suggested mechanism is most visible for the hosts, Russia and Qatar.
%, for which the greatest biases in matchup probabilities substantially decrease.
%Labelling the groups after the draw also makes the Skip mechanism less unfair for Serbia, Poland (the two European teams in Pot 3), and Ecuador.
The result is intuitive: by randomising the group of the host and calling it Group A only after the draw is finished, any bias caused by the sensitivity of the Skip mechanism to the labelling of the groups can be eliminated.

Consequently, we have a clear suggestion for improving the fairness of the FIFA World Cup draw by not assigning the host to Group A automatically.
%Since the groups can be labelled even after he draw is finished, Qatar could be guaranteed to play in Group A without the unfair pre-assignment.

\subsection{The impact of draw constraints} \label{Sec42}

\begin{figure}[t!]
\centering

\begin{subfigure}{\textwidth}
\centering
\caption{2018 FIFA World Cup draw}
\label{Fig2a}

\begin{tikzpicture}
\begin{axis}[
xlabel = {Average number of intra-continental group matches},
x label style = {font=\small},
ylabel = {Proportion of valid assignments (\%, log scale)},
y label style = {font=\small},
ymode = log,
log ticks with fixed point,
y tick label style = {/pgf/number format/.cd,fixed,fixed zerofill,precision=2},
width = \textwidth,
height = 0.6\textwidth,
ymajorgrids = true,
xmin = 5.9,
xmax = 11.1,
ymin = 0.1,
ymax = 110,
%max space between ticks=50,
%legend style = {font=\small,at={(0,-0.16)},anchor=north west,legend columns=1},
%legend entries = {Ratio of feasible solutions without Constraint E$\qquad \qquad \qquad \quad \; \;$, Ratio of feasible solutions with Constraint E$\qquad \qquad \qquad \qquad \quad$}
] 
% Without Constraint E
\addplot [blue, mark=asterisk, only marks, mark size=2pt, mark options={solid,semithick}] coordinates {
(10.62319,100)
(9.276329,35.7137296214757)
(10.297641,65.6189818430356)
(8.941048,23.4349200138148)
(9.701716,35.7117541179784)
(8.352729,12.7510704230809)
(9.396466,24.1014131482362)
(8.03479,8.60464605446038)
(9.965695,49.9984027384196)
(8.617325,17.855055733684)
(9.644732,33.4781525400062)
(8.285341,11.9566222165617)
(9.275426,21.4262003920982)
(7.929512,7.65095460447965)
(8.969907,14.461461227508)
(7.61353,5.16400875196835)
};
% With Constraint E
\addplot [red, mark=pentagon, only marks, mark size=2pt, mark options={solid,semithick}] coordinates {
(8.146898,6.06098120547718)
(7.353034,3.57275890228477)
(7.907861,4.21489758002427)
(7.12505,2.50410167078502)
(7.319901,2.3146404490078)
(6.564022,1.40670361144332)
(7.132526,1.68226689687077)
(6.378716,1.02222023978858)
(7.540045,3.18102316643124)
(6.780325,1.91767231517604)
(7.344857,2.28197887352468)
(6.587272,1.37811655051943)
(6.920189,1.38926743099776)
(6.16852,0.844354525722584)
(6.750677,1.01438677249791)
(6,0.61753782346608)
};
%\addplot [black, mark=asterisk, only marks, mark size=4pt, mark options={solid,semithick}] coordinates {
%(9.276329,35.7137296214757)
%(8.941048,23.4349200138148)
%(8.352729,12.7510704230809)
%(8.03479,8.60464605446038)
%(8.617325,17.855055733684)
%};
\end{axis}
\end{tikzpicture}
\end{subfigure}

\vspace{0.25cm}
\begin{subfigure}{\textwidth}
\centering
\caption{2022 FIFA World Cup draw}
\label{Fig2b}

\begin{tikzpicture}
\begin{axis}[
xlabel = {Average number of intra-continental group matches},
x label style = {font=\small},
ylabel = {Proportion of valid assignments (\%, log scale)},
y label style = {font=\small},
ymode = log,
log ticks with fixed point,
y tick label style = {/pgf/number format/.cd,fixed,fixed zerofill,precision=2},
width = \textwidth,
height = 0.6\textwidth,
ymajorgrids = true,
xmin = 4.9,
xmax = 10.55,
ymin = 0.1,
ymax = 110,
%max space between ticks=50,
legend style = {font=\small,at={(0.2,-0.15)},anchor=north west,legend columns=2},
legend entries = {Without Constraint E$\qquad$, With Constraint E}
] 
% Without Constraint E
\addplot [blue, mark=asterisk, only marks, mark size=2pt, mark options={solid,semithick}] coordinates {
(10.090800568288,100)
(9.02474509887608,26.7857456039543)
(9.57426732475483,53.5682500513694)
(8.52039021169614,14.8785119336672)
(9.1998636410059,35.714981841872)
(8.09616610799281,9.56748302828943)
(8.65002096609122,19.1300839742908)
(7.53675012125127,5.3143572653996)
(8.77780368268262,13.393749056327)
(7.70248998671161,3.45386927516075)
(8.23795141871821,7.17534409706611)
(7.171469,1.92258512140105)
(8.11792514927082,6.1617399486978)
(7.03169959402378,1.61304538895889)
(7.54714969815475,3.30068726246146)
(6.479488,0.900489416700293)
};
% With Constraint E
\addplot [red, mark=pentagon, only marks, mark size=2pt, mark options={solid,semithick}] coordinates {
(7.6440014647344,11.0269107622045)
(6.83541359053917,3.72084601915839)
(7.21717277940177,6.37504735617944)
(7.47718826122263,2.30562407478043)
(6.73669635115813,4.08674395593224)
(6.02597196809429,1.46887946622857)
(6.37528867564028,2.4552796911618)
(5.7359080159006,0.932958180075244)
(6.59334170687361,1.75926706346739)
(5.73879271440304,0.539691852222675)
(6.23617692950027,1.07954854715772)
(5.449065,0.355654996947841)
(6.01958567678563,0.881321716144135)
(5.25366835851235,0.28603511201533)
(5.69713594061169,0.54160451527556)
(5,0.187534371299572)
};
%\addplot [black, mark=asterisk, only marks, mark size=4pt, mark options={solid,semithick}] coordinates {
%(7.47718826122263,2.30562407478043)
%};
\end{axis}
\end{tikzpicture}
\end{subfigure}

%\captionsetup{justification=centerfirst}
\caption{The average number of intra-continental games and the chance of feasibility}
\label{Fig2}
\end{figure}

%\end{document}

Figure~\ref{Fig2} shows the number of intra-continental group matches and the proportion of valid assignments under the 32 sets of draw constraints. Unsurprisingly, there is a strong negative relationship between the two variables. In the case of the 2018 FIFA World Cup, essentially two parallel lines can be seen (but note the logarithmic scale on the vertical axis, the ratio of feasible assignments), where the top line corresponds to the 16 sets of constraints that include Constraint D, showing a more favourable trade-off.

The most restrictive constraint is Constraint E, affecting the UEFA teams, which reduces the probability of feasibility by at least 88\% (2018) and 79\% (2022), respectively. For the 2018 FIFA World Cup, Constraint E in itself is almost as restrictive as Constraints A--D together. A detailed analysis of the restrictiveness of Constraint E in the 2018 FIFA World Cup is provided in the Appendix. 
On the other hand, Constraint C for the CONCACAF nations decreases the number of possible assignments only between 27\% and 44\% (2018), or between 34\% and 46\% (2022).
When all restrictions are imposed, only one out of 162 (532) unconstrained draws can be accepted on average for the 2018 (2022) FIFA World Cup.

In the absence of draw constraints, the expected number of intra-continental group matches is 10.6 and 10.1, respectively. Due to having 14 (13) UEFA teams for eight groups in the 2018 (2022) FIFA World Cup, the lower bounds for the number of these unattractive games are six and five, respectively, even by imposing all Constraints A--E. Constraint E is quite effective in avoiding intra-continental matches as it decreases their expected number by almost 25\% in both cases. However, this does not mean that a combination of Constraints A--D could not reach the same impact with a smaller distortion compared to a uniform draw, which will be investigated in Section~\ref{Sec43}.

\begin{figure}[t!]
\centering

\begin{subfigure}{\textwidth}
\caption{Measure of non-uniformity: Mean absolute bias $\Delta (\mathcal{C})$}
\label{Fig3a}

\begin{tikzpicture}
\begin{axis}[
xlabel = {Proportion of valid assignments (\%, log scale)},
x label style = {font=\small},
ylabel = {Mean absolute bias (percentage point)},
y label style = {font=\small},
%y tick label style = {/pgf/number format/.cd,fixed,fixed zerofill,precision=2},
width = \textwidth,
height = 0.6\textwidth,
ymajorgrids = true,
xmode = log,
log ticks with fixed point,
xmin = 0.5,
xmax = 115,
ymin = 0,
ymax = 1.5,
%max space between ticks=50,
%legend style = {font=\small,at={(0,-0.25)},anchor=north west,legend columns=5},
%legend entries = {15 states$\quad$,17 states$\quad$,20 states}
] 
% Mean distortion without Constraint E
\addplot [blue, mark=asterisk, only marks, mark size=2pt, mark options={solid,semithick}] coordinates {
(100,0.0255573579482397)
(35.7137296214757,0.0264189417500631)
(65.6189818430356,0.0265163304408966)
(23.4349200138148,0.0319429468744676)
(35.7117541179784,0.0287686587378956)
(12.7510704230809,0.042015548142377)
(24.1014131482362,0.0853500876408102)
(8.60464605446038,0.10281751211994)
(49.9984027384196,0.0284658277316085)
(17.855055733684,0.0360599748242196)
(33.4781525400062,0.0961150298294434)
(11.9566222165617,0.111431701176822)
(21.4262003920982,0.0948657188548219)
(7.65095460447965,0.102668809281652)
(14.461461227508,0.092908114034227)
(5.16400875196835,0.104889705925574)
};
% Mean distortion with Constraint E
\addplot [red, mark=pentagon, only marks, mark size=2pt, mark options={solid,semithick}] coordinates {
(6.06098120547718,1.00313507584015)
(3.57275890228477,0.617199565381414)
(4.21489758002427,0.989893048140216)
(2.50410167078502,0.586336684485855)
(2.3146404490078,1.2368747728402)
(1.40670361144332,0.757917734594567)
(1.68226689687077,1.18096790162938)
(1.02222023978858,0.688588842447729)
(3.18102316643124,1.16622562233977)
(1.91767231517604,0.771389712463711)
(2.28197887352468,1.11522178763991)
(1.37811655051943,0.71522444884264)
(1.38926743099776,1.17035760078303)
(0.844354525722584,0.713654538200803)
(1.01438677249791,1.10857612224394)
(0.61753782346608,0.636276712328767)
};
% Mean distortion with Constraint E, host pre-assigned
\addplot [ForestGreen, mark=oplus, only marks, mark size=2pt, mark options={solid,semithick}] coordinates {
(0.61753782346608,0.821769151)
};
\end{axis}
\end{tikzpicture}
\end{subfigure}

\vspace{0.25cm}
\begin{subfigure}{\textwidth}
\caption{Measure of non-uniformity: Maximum absolute bias $\Omega (\mathcal{C})$}
\label{Fig3b}

\begin{tikzpicture}
\begin{axis}[
xlabel = {Proportion of valid assignments (\%, log scale)},
x label style = {font=\small},
ylabel = {Maximum absolute bias (percentage point)},
y label style = {font=\small},
%y tick label style = {/pgf/number format/.cd,fixed,fixed zerofill,precision=2},
width = \textwidth,
height = 0.6\textwidth,
ymajorgrids = true,
xmode = log,
log ticks with fixed point,
xmin = 0.5,
xmax = 115,
ymin = 0,
ymax = 14.5,
%max space between ticks=50,
legend style = {font=\small,at={(0.05,-0.15)},anchor=north west,legend columns=1},
legend entries = {Skip mechanism without Constraint E$\qquad \,\, \qquad \qquad \qquad \qquad \qquad \qquad \quad \; \;$, Skip mechanism with Constraint E$\qquad \,\, \qquad \qquad \qquad \qquad \qquad \qquad \qquad \quad$, Skip mechanism with all constraints if the host is pre-assigned to Group A}
] 
% Max distortion without Constraint E
\addplot [blue, mark=asterisk, only marks, mark size=2pt, mark options={solid,semithick}] coordinates {
(100,0.0940727087625354)
(35.7137296214757,0.0962507705474111)
(65.6189818430356,0.116967008946223)
(23.4349200138148,0.198836053567525)
(35.7117541179784,0.122547165218685)
(12.7510704230809,0.302973783298358)
(24.1014131482362,0.852897811492286)
(8.60464605446038,0.967007678481871)
(49.9984027384196,0.100477193209327)
(17.855055733684,0.154002720713398)
(33.4781525400062,3.29347425512814)
(11.9566222165617,3.4212818625523)
(21.4262003920982,2.91945009709283)
(7.65095460447965,2.9384748683618)
(14.461461227508,3.15493030600135)
(5.16400875196835,3.1252836528544)
};
% Max distortion with Constraint E
\addplot [red, mark=pentagon, only marks, mark size=2pt, mark options={solid,semithick}] coordinates {
(6.06098120547718,10.2517103830224)
(3.57275890228477,12.8715878706903)
(4.21489758002427,10.1050514987487)
(2.50410167078502,10.8210268328772)
(2.3146404490078,10.5560369839623)
(1.40670361144332,13.0420624602763)
(1.68226689687077,10.3721185610788)
(1.02222023978858,10.7296889302103)
(3.18102316643124,13.1111109342807)
(1.91767231517604,13.1475166953054)
(2.28197887352468,11.4641513426566)
(1.37811655051943,11.0489056312625)
(1.38926743099776,10.3350568821988)
(0.844354525722584,12.8246049642651)
(1.01438677249791,10.1159317313505)
(0.61753782346608,10.825)
};
% Max distortion with Constraint E, host pre-assigned
\addplot [ForestGreen, mark=oplus, only marks, mark size=2pt, mark options={solid,semithick}] coordinates {
(0.61753782346608,10.23657)
};
\end{axis}
\end{tikzpicture}
\end{subfigure}

\captionsetup{justification=centering}
\caption{The aggregated bias of the Skip mechanism for the \\ 32 sets of draw constraints, 2018 FIFA World Cup draw}
\label{Fig3}
\end{figure}

%\end{document}

\begin{figure}[t!]
\centering

\begin{subfigure}{\textwidth}
\caption{Measure of non-uniformity: Mean absolute bias $\Delta (\mathcal{C})$}
\label{Fig4a}

\begin{tikzpicture}
\begin{axis}[
xlabel = {Proportion of valid assignments (\%, log scale)},
x label style = {font=\small},
ylabel = {Mean absolute bias (percentage point)},
y label style = {font=\small},
%y tick label style = {/pgf/number format/.cd,fixed,fixed zerofill,precision=2},
width = \textwidth,
height = 0.6\textwidth,
ymajorgrids = true,
xmode = log,
log ticks with fixed point,
xmin = 0.09,
xmax = 115,
ymin = 0,
ymax = 1.5,
%max space between ticks=50,
%legend style = {font=\small,at={(0,-0.25)},anchor=north west,legend columns=5},
%legend entries = {15 states$\quad$,17 states$\quad$,20 states}
] 
% Mean distortion without Constraint E
\addplot [blue, mark=asterisk, only marks, mark size=2pt, mark options={solid,semithick}] coordinates {
(100,0.0282456071846973)
(26.7857456039543,0.0236182486686763)
(53.5682500513694,0.0267210119017049)
(14.8785119336672,0.206564501278637)
(35.714981841872,0.0273673815416241)
(9.56748302828943,0.314594456557567)
(19.1300839742908,0.227458884007567)
(5.3143572653996,0.544662042991035)
(13.393749056327,0.0287994958527391)
(3.45386927516075,0.300263355644696)
(7.17534409706611,0.264066761034612)
(1.92258512140105,0.553531455260589)
(6.1617399486978,0.201475762000948)
(1.61304538895889,0.569971289164597)
(3.30068726246146,0.402648567605061)
(0.900489416700293,0.711261590125836)
};
% Mean distortion with Constraint E
\addplot [red, mark=pentagon, only marks, mark size=2pt, mark options={solid,semithick}] coordinates {
(11.0269107622045,0.976679369130177)
(3.72084601915839,0.983432534866167)
(6.37504735617944,0.973863789860189)
(2.30562407478043,1.20130754078856)
(4.08674395593224,1.16866239586788)
(1.46887946622857,1.05382261646784)
(2.4552796911618,1.1075922909458)
(0.932958180075244,0.99023771745208)
(1.75926706346739,1.14466281205177)
(0.539691852222675,1.36640647256331)
(1.07954854715772,1.08359799422733)
(0.355654996947841,1.38927318479867)
(0.881321716144135,1.09707546907179)
(0.28603511201533,1.27211678049165)
(0.54160451527556,1.0313958253772)
(0.187534371299572,1.26630591549296)
};
% Mean distortion with Constraint E, host pre-assigned
\addplot [ForestGreen, mark=oplus, only marks, mark size=2pt, mark options={solid,semithick}] coordinates {
(0.187534371299572,1.433537746)
};
\end{axis}
\end{tikzpicture}
\end{subfigure}

\vspace{0.25cm}
\begin{subfigure}{\textwidth}
\caption{Measure of non-uniformity: Maximum absolute bias $\Omega (\mathcal{C})$}
\label{Fig4b}

\begin{tikzpicture}
\begin{axis}[
xlabel = {Proportion of valid assignments (\%, log scale)},
x label style = {font=\small},
ylabel = {Maximum absolute bias (percentage point)},
y label style = {font=\small},
%y tick label style = {/pgf/number format/.cd,fixed,fixed zerofill,precision=2},
width = \textwidth,
height = 0.6\textwidth,
ymajorgrids = true,
xmode = log,
log ticks with fixed point,
xmin = 0.09,
xmax = 115,
ymin = 0,
ymax = 14.5,
%max space between ticks=50,
legend style = {font=\small,at={(0.05,-0.15)},anchor=north west,legend columns=1},
legend entries = {Skip mechanism without Constraint E$\qquad \,\, \qquad \qquad \qquad \qquad \qquad \qquad \quad \; \;$, Skip mechanism with Constraint E$\qquad \,\, \qquad \qquad \qquad \qquad \qquad \qquad \qquad \quad$, Skip mechanism with all constraints if the host is pre-assigned to Group A}
] 
% Max distortion without Constraint E
\addplot [blue, mark=asterisk, only marks, mark size=2pt, mark options={solid,semithick}] coordinates {
(100,0.105072629871603)
(26.7857456039543,0.0964104804641083)
(53.5682500513694,0.123464298328568)
(14.8785119336672,1.77925840787562)
(35.714981841872,0.107831775550182)
(9.56748302828943,1.43162828530223)
(19.1300839742908,1.78416927418194)
(5.3143572653996,2.25122239424852)
(13.393749056327,0.109104988490938)
(3.45386927516075,2.65090251600944)
(7.17534409706611,2.36692069619457)
(1.92258512140105,3.44081937153552)
(6.1617399486978,3.18815439426307)
(1.61304538895889,2.88790701974996)
(3.30068726246146,3.19788622594584)
(0.900489416700293,3.08211990553404)
};
% Max distortion with Constraint E
\addplot [red, mark=pentagon, only marks, mark size=2pt, mark options={solid,semithick}] coordinates {
(11.0269107622045,6.25032240478187)
(3.72084601915839,11.8512693813907)
(6.37504735617944,9.19892525122285)
(2.30562407478043,11.3920508860288)
(4.08674395593224,6.37778438896285)
(1.46887946622857,11.2676226968149)
(2.4552796911618,8.56170325455368)
(0.932958180075244,11.2067461907007)
(1.75926706346739,7.87908560813846)
(0.539691852222675,11.2588126541223)
(1.07954854715772,8.29407469204141)
(0.355654996947841,9.26444601458809)
(0.881321716144135,8.13352135943216)
(0.28603511201533,10.607584358603)
(0.54160451527556,8.42631006818493)
(0.187534371299572,8.5581)
};
% Max distortion with Constraint E, host pre-assigned
\addplot [ForestGreen, mark=oplus, only marks, mark size=2pt, mark options={solid,semithick}] coordinates {
(0.187534371299572,10.49049)
};
\end{axis}
\end{tikzpicture}
\end{subfigure}

\captionsetup{justification=centering}
\caption{The aggregated bias of the Skip mechanism for the \\ 32 sets of draw constraints, 2022 FIFA World Cup draw}
\label{Fig4}
\end{figure}

%\end{document}

Figures~\ref{Fig3} and \ref{Fig4} present the distortion of the draw under all sets of draw constraints as a function of the proportion of valid assignments.
% : Figure~\ref{Fig3a} if the distortion is measured by mean absolute bias $\Delta (\mathcal{C})$ and Figure~\ref{Fig3b} if it is measured by maximum absolute bias $\Omega (\mathcal{C})$
Imposing more constraints does not necessarily imply a higher bias. The draw procedure clearly struggles with Constraint E, especially in the 2018 FIFA World Cup draw. If the number of UEFA teams is not restricted in any group, $\Delta (\mathcal{C})$ ($\Omega (\mathcal{C})$) never exceeds 0.12 and 0.72 (3.5 and 3.5) for the 2018 and 2022 draws, respectively, but $\Delta (\mathcal{C})$ ($\Omega (\mathcal{C})$) is certainly higher than 0.58 and 0.97 (10.1 and 6.2), respectively, if Constraint E is imposed. Constraint E is more challenging for the Skip mechanism since it allows more than one team from a given set to be assigned to the same group \citep{Csato2026a}. On the other hand, Constraints A--D can be guaranteed by prohibited pairs, and imposing only one of them preserves uniformity for the Skip mechanism \citep[Lemma~1]{Csato2026a}. These results are insensitive to the measure of non-uniformity; the patterns seen in Figures~\ref{Fig3a} and \ref{Fig3b}, as well as Figures~\ref{Fig4a} and \ref{Fig4b} are similar.

\subsection{The Pareto frontier of the draw constraints} \label{Sec43}

\begin{figure}[t!]
\centering

\begin{subfigure}{\textwidth}
\caption{Measure of non-uniformity: Mean absolute bias $\Delta (\mathcal{C})$}
\label{Fig5a}

\begin{tikzpicture}
\begin{axis}[
name = axis1,
xlabel = Average number of intra-continental group matches,
x label style = {font=\small},
ylabel = Mean absolute bias (percentage point),
y label style = {font=\small},
yticklabel style = {/pgf/number format/fixed,/pgf/number format/precision=5},
scaled y ticks = false,
width = \textwidth,
height = 0.6\textwidth,
ymajorgrids = true,
xmin = 5.9,
xmax = 11.1,
ymin = 0,
%ymax = 0.055,
%max space between ticks=50,
%legend style = {font=\small,at={(0,-0.15)},anchor=north west,legend columns=1},
%legend entries = {Uniform mechanism without Constraints A and E$\qquad \qquad \qquad \quad$, Uniform mechanism with Constraint A and without Constraint E, Uniform mechanism with Constraint E and without Constraint A, Uniform mechanism with Constraints A and E$\qquad \qquad \qquad \qquad \;$, Skip mechanism without Constraints A and E$\qquad \qquad \qquad \qquad \; \,$, Skip mechanism with Constraint A and without Constraint E$\quad \; \,$, Skip mechanism with Constraint E$\qquad \qquad \qquad \qquad \qquad \quad \qquad \;$}
] 
% Skip mechanism without Constraint E
\addplot [blue, mark=asterisk, only marks, mark size=2pt, mark options={solid,thin}] coordinates {
(10.62319,0)
(9.276329,0)
(10.297641,0)
(8.941048,0.0319429468744676)
(9.701716,0)
(8.352729,0.042015548142377)
(9.396466,0.0853500876408102)
(8.03479,0.10281751211994)
(9.965695,0)
(8.617325,0.0360599748242196)
(9.644732,0.0961150298294434)
(8.285341,0.111431701176822)
(9.275426,0.0948657188548219)
(7.929512,0.102668809281652)
(8.969907,0.092908114034227)
(7.61353,0.104889705925574)
};
% Skip mechanism with Constraint E
\addplot [red, mark=pentagon, only marks, mark size=2pt, mark options={solid,semithick}] coordinates {
(8.146898,1.00313507584015)
(7.353034,0.617199565381414)
(7.907861,0.989893048140216)
(7.12505,0.586336684485855)
(7.319901,1.2368747728402)
(6.564022,0.757917734594567)
(7.132526,1.18096790162938)
(6.378716,0.688588842447729)
(7.540045,1.16622562233977)
(6.780325,0.771389712463711)
(7.344857,1.11522178763991)
(6.587272,0.71522444884264)
(6.920189,1.17035760078303)
(6.16852,0.713654538200803)
(6.750677,1.10857612224394)
(6,0.636276712328767)
};
% Current voting power
\draw [black,thick,dotted] (6,0.636276712328767) -- (7.61353,0.104889705925574) -- (8.352729,0.042015548142377) -- (9.276329,0);
% Optimal line
%\draw [brown,thick,dotted] (5,1.26630591549296) -- (8.77780368268262,0);
\end{axis}
\end{tikzpicture}
\end{subfigure}

\vspace{0.25cm}
\begin{subfigure}{\textwidth}
\caption{Measure of non-uniformity: Maximum absolute bias $\Omega (\mathcal{C})$}
\label{Fig5b}

\begin{tikzpicture}
\begin{axis}[
name = axis1,
xlabel = Average number of intra-continental group matches,
x label style = {font=\small},
ylabel = {Maximum absolute bias (percentage point)},
y label style = {font=\small},
yticklabel style = {/pgf/number format/fixed,/pgf/number format/precision=5},
scaled y ticks = false,
width = \textwidth,
height = 0.6\textwidth,
ymajorgrids = true,
xmin = 5.9,
xmax = 11.1,
ymin = 0,
%ymax = 0.055,
%max space between ticks=50,
legend style = {font=\small,at={(0,-0.15)},anchor=north west,legend columns=2},
legend entries = {Skip mechanism without Constraint E$\qquad$, Skip mechanism with Constraint E}
] 
% Skip mechanism without Constraint E
\addplot [blue, mark=asterisk, only marks, mark size=2pt, mark options={solid,thin}] coordinates {
(10.62319,0)
(9.276329,0)
(10.297641,0)
(8.941048,0.198836053567525)
(9.701716,0)
(8.352729,0.302973783298358)
(9.396466,0.852897811492286)
(8.03479,0.967007678481871)
(9.965695,0)
(8.617325,0.154002720713398)
(9.644732,3.29347425512814)
(8.285341,3.4212818625523)
(9.275426,2.91945009709283)
(7.929512,2.9384748683618)
(8.969907,3.15493030600135)
(7.61353,3.1252836528544)
};
% Skip mechanism with Constraint E
\addplot [red, mark=pentagon, only marks, mark size=2pt, mark options={solid,semithick}] coordinates {
(8.146898,10.2517103830224)
(7.353034,12.8715878706903)
(7.907861,10.1050514987487)
(7.12505,10.8210268328772)
(7.319901,10.5560369839623)
(6.564022,13.0420624602763)
(7.132526,10.3721185610788)
(6.378716,10.7296889302103)
(7.540045,13.1111109342807)
(6.780325,13.1475166953054)
(7.344857,11.4641513426566)
(6.587272,11.0489056312625)
(6.920189,10.3350568821988)
(6.16852,12.8246049642651)
(6.750677,10.1159317313505)
(6,10.825)
};
% Pareto Frontier
\draw [black,thick,dotted] (6,10.825) -- (8.03479,0.967007678481871) -- (8.352729,0.302973783298358) -- (8.617325,0.154002720713398) -- (9.276329,0)
;
\end{axis}
\end{tikzpicture}
\end{subfigure}

\captionsetup{justification=centering}
\caption{The trade-off between intra-continental group \\ matches and non-uniformity, 2018 FIFA World Cup \vspace{0.2cm} \\
\footnotesize{\emph{Note}: The dotted black line shows the Pareto frontier.}}
\label{Fig5} 
\end{figure}

%\end{document}

\begin{figure}[t!]
\centering

\begin{subfigure}{\textwidth}
\caption{Measure of non-uniformity: Mean absolute bias $\Delta (\mathcal{C})$}
\label{Fig6a}

\begin{tikzpicture}
\begin{axis}[
name = axis1,
xlabel = Average number of intra-continental group matches,
x label style = {font=\small},
ylabel = Mean absolute bias (percentage point),
y label style = {font=\small},
yticklabel style = {/pgf/number format/fixed,/pgf/number format/precision=5},
scaled y ticks = false,
width = \textwidth,
height = 0.6\textwidth,
ymajorgrids = true,
xmin = 4.9,
xmax = 10.55,
ymin = 0,
%ymax = 0.055,
%max space between ticks=50,
%legend style = {font=\small,at={(0,-0.15)},anchor=north west,legend columns=1},
%legend entries = {Uniform mechanism without Constraints A and E$\qquad \qquad \qquad \quad$, Uniform mechanism with Constraint A and without Constraint E, Uniform mechanism with Constraint E and without Constraint A, Uniform mechanism with Constraints A and E$\qquad \qquad \qquad \qquad \;$, Skip mechanism without Constraints A and E$\qquad \qquad \qquad \qquad \; \,$, Skip mechanism with Constraint A and without Constraint E$\quad \; \,$, Skip mechanism with Constraint E$\qquad \qquad \qquad \qquad \qquad \quad \qquad \;$}
] 
% Skip mechanism without Constraint E
\addplot [blue, mark=asterisk, only marks, mark size=2pt, mark options={solid,thin}] coordinates {
(10.090800568288,0)
(9.02474509887608,0)
(9.57426732475483,0)
(8.52039021169614,0.206564501278637)
(9.1998636410059,0)
(8.09616610799281,0.314594456557567)
(8.65002096609122,0.227458884007567)
(7.53675012125127,0.544662042991035)
(8.77780368268262,0)
(7.70248998671161,0.300263355644696)
(8.23795141871821,0.264066761034612)
(7.171469,0.553531455260589)
(8.11792514927082,0.201475762000948)
(7.03169959402378,0.569971289164597)
(7.54714969815475,0.402648567605061)
(6.479488,0.711261590125836)
};
% Skip mechanism with Constraint E
\addplot [red, mark=pentagon, only marks, mark size=2pt, mark options={solid,semithick}] coordinates {
(7.6440014647344,0.976679369130177)
(6.83541359053917,0.983432534866167)
(7.21717277940177,0.973863789860189)
(7.47718826122263,1.20130754078856)
(6.73669635115813,1.16866239586788)
(6.02597196809429,1.05382261646784)
(6.37528867564028,1.1075922909458)
(5.7359080159006,0.99023771745208)
(6.59334170687361,1.14466281205177)
(5.73879271440304,1.36640647256331)
(6.23617692950027,1.08359799422733)
(5.449065,1.38927318479867)
(6.01958567678563,1.09707546907179)
(5.25366835851235,1.27211678049165)
(5.69713594061169,1.0313958253772)
(5,1.26630591549296)
};
% Current voting power
\draw [black,thick,dotted] (5,1.26630591549296) -- (6.479488,0.711261590125836) -- (7.70248998671161,0.300263355644696) -- (8.77780368268262,0);
% Optimal line
%\draw [brown,thick,dotted] (5,1.26630591549296) -- (8.77780368268262,0);
\end{axis}
\end{tikzpicture}
\end{subfigure}

\vspace{0.25cm}
\begin{subfigure}{\textwidth}
\caption{Measure of non-uniformity: Maximum absolute bias $\Omega (\mathcal{C})$}
\label{Fig6b}

\begin{tikzpicture}
\begin{axis}[
name = axis1,
xlabel = Average number of intra-continental group matches,
x label style = {font=\small},
ylabel = {Maximum absolute bias (percentage point)},
y label style = {font=\small},
yticklabel style = {/pgf/number format/fixed,/pgf/number format/precision=5},
scaled y ticks = false,
width = 0.99\textwidth,
height = 0.6\textwidth,
ymajorgrids = true,
xmin = 4.9,
xmax = 10.55,
ymin = 0,
%ymax = 0.055,
%max space between ticks=50,
legend style = {font=\small,at={(0,-0.15)},anchor=north west,legend columns=2},
legend entries = {Skip mechanism without Constraint E$\qquad$, Skip mechanism with Constraint E}
] 
% Skip mechanism without Constraint E
\addplot [blue, mark=asterisk, only marks, mark size=2pt, mark options={solid,thin}] coordinates {
(10.090800568288,0)
(9.02474509887608,0)
(9.57426732475483,0)
(8.52039021169614,1.77925840787562)
(9.1998636410059,0)
(8.09616610799281,1.43162828530223)
(8.65002096609122,1.78416927418194)
(7.53675012125127,2.25122239424852)
(8.77780368268262,0)
(7.70248998671161,2.65090251600944)
(8.23795141871821,2.36692069619457)
(7.171469,3.44081937153552)
(8.11792514927082,3.18815439426307)
(7.03169959402378,2.88790701974996)
(7.54714969815475,3.19788622594584)
(6.479488,3.08211990553404)
};
% Skip mechanism with Constraint E
\addplot [red, mark=pentagon, only marks, mark size=2pt, mark options={solid,semithick}] coordinates {
(7.6440014647344,6.25032240478187)
(6.83541359053917,11.8512693813907)
(7.21717277940177,9.19892525122285)
(7.47718826122263,11.3920508860288)
(6.73669635115813,6.37778438896285)
(6.02597196809429,11.2676226968149)
(6.37528867564028,8.56170325455368)
(5.7359080159006,11.2067461907007)
(6.59334170687361,7.87908560813846)
(5.73879271440304,11.2588126541223)
(6.23617692950027,8.29407469204141)
(5.449065,9.26444601458809)
(6.01958567678563,8.13352135943216)
(5.25366835851235,10.607584358603)
(5.69713594061169,8.42631006818493)
(5,8.5581)
};
% Pareto Frontier
\draw [black,thick,dotted] (5,8.5581) -- (6.479488,3.08211990553404) -- (8.77780368268262,0);
\end{axis}
\end{tikzpicture}
\end{subfigure}

\captionsetup{justification=centering}
\caption{The trade-off between intra-continental group \\ matches and non-uniformity, 2022 FIFA World Cup \vspace{0.2cm} \\
\footnotesize{\emph{Note}: The dotted black line shows the Pareto frontier.}}
\label{Fig6} 
\end{figure}

%\end{document}

Figures~\ref{Fig5} and \ref{Fig6} uncover the trade-off between the number of intra-continental group matches and the non-uniformity of the draw for the 2018 and 2022 FIFA World Cups, respectively. The expected number of unattractive games is 10.6 and 10.1, respectively, in an unconstrained draw that can be decreased to 9.3 and 8.8, respectively, without any price by imposing one of Constraints A--D. If non-uniformity is measured by the average distortion, the Pareto frontier contains the draw under Constraints A--D for both FIFA World Cups, which implies about 7.6 (2018) or 6.5 (2022) intra-continental group matches for a moderate level of non-uniformity. Finally, imposing all restrictions is also efficient as the number of intra-continental games reaches its minimum of six (2018) or five (2022), albeit the price may be excessive.

\begin{table}[t!]
  \centering
  \caption{The optimal (Pareto efficient) sets of draw constraints}  
  \label{Table4}
  
\begin{subtable}{\textwidth}
\centering
\caption{2018 FIFA World Cup draw}
\label{Table4a}
  \rowcolors{3}{}{gray!20}
    \begin{tabularx}{\textwidth}{CcCc} \toprule
    Problem & Measure of non-uniformity & Optimal $\mathcal{C}$ & Value of $\alpha$ \\ \bottomrule
    \eqref{Opt1} & $\Delta(\mathcal{C})$ & 2 (D) & $0.6873 \leq \alpha$ \\
    \eqref{Opt1} & $\Delta(\mathcal{C})$ & 10 (B, D) & $0.5404 \leq \alpha \leq 0.6873$ \\
    \eqref{Opt1} & $\Delta(\mathcal{C})$ & 30 (A--D) & $0.2329 \leq \alpha \leq 0.5404$ \\
    \eqref{Opt1} & $\Delta(\mathcal{C})$ & 31 (A--E) & $\alpha \leq 0.2329$ \\ \hline
    \eqref{Opt2} & $\Omega(\mathcal{C})$ & 2 (D) & $0.8106 \leq \alpha$ \\
    \eqref{Opt2} & $\Omega(\mathcal{C})$ & 18 (A, D) & $0.6398 \leq \alpha \leq 0.8106$ \\
    \eqref{Opt2} & $\Omega(\mathcal{C})$ & 10 (B, D) & $0.3238 \leq \alpha \leq 0.6398$ \\
    \eqref{Opt2} & $\Omega(\mathcal{C})$ & 14 (B, C, D) & $0.1711 \leq \alpha \leq 0.3238$ \\
    \eqref{Opt2} & $\Omega(\mathcal{C})$ & 31 (A--E) & $\alpha \leq 0.1711$ \\ \toprule
    \end{tabularx}
\end{subtable}

\vspace{0.25cm}
\begin{subtable}{\textwidth}
\centering
\caption{2022 FIFA World Cup draw}
\label{Table4b}
  \rowcolors{3}{}{gray!20}
    \begin{tabularx}{\textwidth}{CcCc} \toprule
    Problem & Measure of non-uniformity & Optimal $\mathcal{C}$ & Value of $\alpha$ \\ \bottomrule
    \eqref{Opt1} & $\Delta(\mathcal{C})$ & 16 (A) & $0.2637 \leq \alpha$ \\
    \eqref{Opt1} & $\Delta(\mathcal{C})$ & 18 (A, D) & $0.2293 \leq \alpha \leq 0.2637$ \\
    \eqref{Opt1} & $\Delta(\mathcal{C})$ & 30 (A--D) & $0.2105 \leq \alpha \leq 0.2293$ \\
    \eqref{Opt1} & $\Delta(\mathcal{C})$ & 31 (A--E) & $\alpha \leq 0.2105$ \\ \hline
    \eqref{Opt2} & $\Omega(\mathcal{C})$ & 16 (A) & $0.4272 \leq \alpha$ \\
    \eqref{Opt2} & $\Omega(\mathcal{C})$ & 30 (A--D) & $0.2127 \leq \alpha \leq 0.4272$ \\
    \eqref{Opt2} & $\Omega(\mathcal{C})$ & 31 (A--E) & $\alpha \leq 0.2127$ \\ \toprule
    \end{tabularx}
\end{subtable}

\end{table}

Table~\ref{Table4} summarises the optimal (Pareto efficient) solutions to problems~\eqref{Opt1} and \eqref{Opt2} for the 2018 (Table~\ref{Table4a}) and 2022 (Table~\ref{Table4b}) FIFA World Cups as a function of parameter $\alpha$. Besides the three sets of constraints discussed above, at most three further sets can be efficient. In addition, one scenario is an optimal solution to problem~\eqref{Opt1} only in a narrow range of parameter $\alpha$ for both the 2018 (scenario 14 with Constraints A and D) and 2022 (scenario 30 with Constraints A--D) FIFA World Cup draws.

\input{Figure7_objective_value_alpha_2018}

\input{Figure8_objective_value_alpha_2022}

However, reporting only the efficient set of draw constraints for a given value of parameter $\alpha$ may be misleading because another marginally inefficient set could be overlooked. Therefore, Figures~\ref{Fig7} and \ref{Fig8} show the values of the objective function in optimisation problems~\eqref{Opt1} and \eqref{Opt1} for all sets of draw constraints. The blue lines (no constraint for UEFA teams) are generally preferred if $\alpha$ is high, but the red lines (Constraint E is used) are better if $\alpha$ is low. Interestingly, focusing on the mean absolute bias $\Delta (\mathcal{C})$ in the 2018 FIFA World Cup draw implies that Constraint E should not be used without Constraint D, as the red lines clearly outperform the black lines.
In both the 2018 and 2022 FIFA World Cup draws, the Pareto frontier is well approximated by only two sets of constraints: all constraints A--E (the thick brown line) and only one particular constraint, Constraint D/A for the 2018/2022 FIFA World Cup (the thick green line).

According to Section~\ref{Sec1}, the current policies of FIBA in basketball and FIFA in football indicate that they prefer a high number of intra-confederation group games to non-uniformity with a small value of $\alpha$. On the other hand, the choices of IHF in handball and FIVB in volleyball are inefficient in the proposed framework even if $\alpha = 1$ since the maximal number of teams from one particular continent can be one in a group without compromising uniformity as shown by \citet[Lemma~1]{Csato2026a}. Figures~\ref{Fig7} and \ref{Fig8} uncover that compromise solutions with some sets of geographic constraints could not result in a substantial (if any) improvement, at least in the case of the 2018 and 2022 FIFA World Cup draws.

Therefore, while the official rule of the FIFA World Cup draw of imposing all geographic constraints can obviously be rationalised if the number intra-continental group games should be minimised, the results of our analysis remain important for the organiser: further sets of draw restrictions might also be chosen if the non-uniformity of the draw is judged excessive and some stakeholders, especially the national teams, call for less unfair rules.

\section{Concluding remarks} \label{Sec6}

Tournament design can cause severe headaches for organisers due to the conflicts between commercial goals and the fairness of the rules. First in the literature, our paper analyses such a trade-off between attractiveness and non-uniformity in the relatively unexplored field of draw procedures. The suggested framework is based on a parametric optimisation model, which allows for determining the Pareto frontier implied by reasonable sets of constraints in a group draw.

The proposed approach is applied to the 2018 and 2022 FIFA World Cup draws with the following findings:
\begin{itemize}
\item
The unnecessary policy of pre-assigning the host to Group A increases the mean absolute bias for this national team by 555\% and 80\%, respectively, while it is essentially not fairer to any other team (Section~\ref{Sec41});
\item
The Skip mechanism struggles with Constraint E concerning UEFA teams, which increases the mean (maximum) absolute bias $\Delta (\mathcal{C})$ ($\Omega (\mathcal{C})$) of the draw by at least 507\% and 78\% (221\% and 155\%) for the 2018 and 2022 tournaments, respectively (Section~\ref{Sec42});
\item
Out of 32 reasonable scenarios, at most five sets of draw constraints provide an optimal trade-off between the expected number of intra-continental group matches and the non-uniformity of the Skip mechanism, and the Pareto frontier can be approximated well by only two extremal sets of draw constraints (Section~\ref{Sec43}).
\end{itemize}

The suboptimality of the draw procedure for the two recent FIFA World Cups has been revealed.
Consequently, tournament organisers are encouraged to revise the pre-assignment of the host(s), as well as to carry out similar analyses before the group draw of major competitions in order to compute the efficient choices with respect to the inevitable trade-off between the attractiveness and uniqueness of group matches and the fairness of the draw.

\section*{Acknowledgements}
\addcontentsline{toc}{section}{Acknowledgements}
\noindent
This paper could not have been written without my \emph{father} (also called \emph{L\'aszl\'o Csat\'o}), who helped to code the simulations in Python. \\
%We are grateful to \emph{Dries Goossens}, \emph{Alex Krumer}, \emph{Frits C.~R.~Spieksma}, and \emph{Stephan Westphal} for useful advice. \\
Two anonymous reviewers provided valuable comments and suggestions on earlier drafts. \\
The research was supported by the National Research, Development and Innovation Office under Grants Advanced 152220 and FK 145838, and the J\'anos Bolyai Research Scholarship of the Hungarian Academy of Sciences.

\bibliographystyle{apalike}
\bibliography{All_references}

\clearpage
\section*{Appendix}
\addcontentsline{toc}{section}{Appendix}

\subsection*{The restrictiveness of Constraint E in the 2018 FIFA World Cup}
\addcontentsline{toc}{subsection}{The restrictiveness of Constraint E in the 2018 FIFA World Cup}

Constraint E is the union of two constraints:
\begin{itemize}
\item
\emph{Constraint E1}: each group contains at least one UEFA team;
\item
\emph{Constraint E2}: each group contains at most two UEFA teams.
\end{itemize}
Here we compute how restrictive Constraints E1, E2, and E are in the 2018 FIFA World Cup draw.

The pots are considered in the order Pot 4, Pot 3, Pot 1, and Pot 2, which does not affect the matchup probabilities in a uniform draw, but simplifies our analysis. According to Table~\ref{Table2a}, Pot 1 contains six, Pot 2 contains four, Pot 3 contains three, and Pot 4 contains one UEFA team(s).

When Pots 3 and 4 with four UEFA teams are emptied, two scenarios can occur:
\begin{itemize}
\item
\emph{One group contains two UEFA teams and two groups contain one UEFA team, which has a probability of 3/8.}

The two non-UEFA teams in Pot 1 can be assigned to the groups in five different ways:
\begin{enumerate}
\item % a
\emph{Both non-UEFA teams in Pot 1 are assigned to two of the five groups without a UEFA team.}

The probability of this event is
\[
\frac{3}{8} \cdot \left( \frac{5}{8} \cdot \frac{4}{7} \right) = \frac{3}{8} \cdot \frac{5}{14} = \frac{15}{112}.
\]
Then, one group has three UEFA teams, two groups have two UEFA teams, and three groups have one UEFA team when Pot 1 is emptied.

Constraint E1 is satisfied if the four non-UEFA teams in Pot 2 are all assigned to four of the six groups with a UEFA team, which has a probability of
\[
\frac{15}{112} \cdot \left( \frac{6}{8} \cdot \frac{5}{7} \cdot \frac{4}{6} \cdot \frac{3}{5} \right) = \frac{15}{112} \cdot \frac{3}{14} = \frac{45}{1568}.
\]

Constraint E2 is certainly violated.

Constraint E (the union of Constraints E1 and E2) is certainly violated.

\item % b
\emph{One non-UEFA team in Pot 1 is assigned to one of the five groups without a UEFA team, and the other non-UEFA team in Pot 1 is assigned to one of the two groups with one UEFA team.}

The probability of this event is
\[
\frac{3}{8} \cdot \left( 2 \cdot \frac{5}{8} \cdot \frac{2}{7} \right) = \frac{3}{8} \cdot \frac{5}{14} = \frac{15}{112}.
\]
Then, one group has three UEFA teams, one group has two UEFA teams, and five groups have one UEFA team when Pot 1 is emptied.

Constraint E1 is satisfied if the four non-UEFA teams in Pot 2 are all assigned to four of the seven groups with a UEFA team, which has a probability of
\[
\frac{15}{112} \cdot \left( \frac{7}{8} \cdot \frac{6}{7} \cdot \frac{5}{6} \cdot \frac{4}{5} \right) = \frac{15}{112} \cdot \frac{1}{2} = \frac{15}{224}.
\]

Constraint E2 is certainly violated.

Constraint E (the union of Constraints E1 and E2) is certainly violated.

\item % c
\emph{One non-UEFA team in Pot 1 is assigned to one of the five groups without a UEFA team, and the other non-UEFA team in Pot 1 is assigned to the only group with two UEFA teams.}

The probability of this event is
\[
\frac{3}{8} \cdot \left( 2 \cdot \frac{5}{8} \cdot \frac{1}{7} \right) = \frac{3}{8} \cdot \frac{5}{28} = \frac{15}{224}.
\]
Then, three groups have two UEFA teams, and four groups have one UEFA team when Pot 1 is emptied.

Constraint E1 is satisfied if the four non-UEFA teams in Pot 2 are all assigned to four of the seven groups with a UEFA team, which has a probability of
\[
\frac{15}{224} \cdot \left( \frac{7}{8} \cdot \frac{6}{7} \cdot \frac{5}{6} \cdot \frac{4}{5} \right) = \frac{15}{224} \cdot \frac{1}{2} = \frac{15}{448}.
\]

Constraint E2 is satisfied if the four UEFA teams in Pot 2 are all assigned to four of the five groups with at most one UEFA team, which has a probability of
\[
\frac{15}{224} \cdot \left( \frac{5}{8} \cdot \frac{4}{7} \cdot \frac{3}{6} \cdot \frac{2}{5} \right) = \frac{15}{224} \cdot \frac{1}{14} = \frac{15}{3136}.
\]

Constraint E (the union of Constraint E1 and E2) is violated if either Constraint E2 is violated, or the four UEFA teams in Pot 2 are all assigned to the four groups with one UEFA team. The latter event has a probability of
\[
\frac{15}{224} \cdot \left( \frac{4}{8} \cdot \frac{3}{7} \cdot \frac{2}{6} \cdot \frac{1}{5} \right) = \frac{15}{224} \cdot \frac{1}{70} = \frac{3}{3136}.
\]
Thus, Constraint E is satisfied with a probability of 
\[
\frac{15}{3136} - \frac{3}{3136} = \frac{12}{3136}.
\]

\item % d
\emph{Both non-UEFA teams in Pot 1 are assigned to the two groups with one UEFA team.}

The probability of this event is
\[
\frac{3}{8} \cdot \left( \frac{2}{8} \cdot \frac{1}{7} \right) = \frac{3}{8} \cdot \frac{1}{28} = \frac{3}{224}.
\]
Then, one group has three UEFA teams, and seven groups have one UEFA team when Pot 1 is emptied.

Constraint E1 is certainly satisfied.

Constraint E2 is certainly violated.

Constraint E (the union of Constraints E1 and E2) is certainly violated.

\item % e
\emph{One non-UEFA team in Pot 1 is assigned to one of the two groups with one UEFA team, and the other non-UEFA team in Pot 1 is assigned to the only group with two UEFA teams.}

The probability of this event is
\[
\frac{3}{8} \cdot \left( 2 \cdot \frac{2}{8} \cdot \frac{1}{7} \right) = \frac{3}{8} \cdot \frac{1}{14} = \frac{3}{112}.
\]
Then, two groups have two UEFA teams, and six groups have one UEFA team when Pot 1 is emptied.

Constraint E1 is certainly satisfied.

Constraint E2 is satisfied if the four UEFA teams in Pot 2 are all assigned to four of the six groups with one UEFA team, which has a probability of
\[
\frac{3}{112} \cdot \left( \frac{6}{8} \cdot \frac{5}{7} \cdot \frac{4}{6} \cdot \frac{3}{5} \right) = \frac{3}{112} \cdot \frac{3}{14} = \frac{9}{1568}.
\]

Constraint E (the union of Constraints E1 and E2) is satisfied if and only if Constraint E2 is satisfied.
\end{enumerate}

\item
\emph{Four groups contain one UEFA team, which has a probability of 5/8.}

The two non-UEFA teams in Pot 1 can be assigned to the groups in three different ways:
\begin{enumerate}
\item % a
\emph{Both non-UEFA teams in Pot 1 are assigned to two of the four groups without a UEFA team.}

The probability of this event is
\[
\frac{5}{8} \cdot \left( \frac{4}{8} \cdot \frac{3}{7} \right) = \frac{5}{8} \cdot \frac{3}{14} = \frac{15}{112}.
\]
Then, four groups have two UEFA teams, and two groups have one UEFA team when Pot 1 is emptied.

Constraint E1 is satisfied if the four non-UEFA teams in Pot 2 are all assigned to four of the six groups with a UEFA team, which has a probability of
\[
\frac{15}{112} \cdot \left( \frac{6}{8} \cdot \frac{5}{7} \cdot \frac{4}{6} \cdot \frac{3}{5} \right) = \frac{15}{112} \cdot \frac{3}{14} = \frac{45}{1568}.
\]

Constraint E2 is satisfied if the four UEFA teams in Pot 2 are all assigned to the four groups with at most one UEFA team, which has a probability of
\[
\frac{15}{112} \cdot \left( \frac{4}{8} \cdot \frac{3}{7} \cdot \frac{2}{6} \cdot \frac{1}{5} \right) = \frac{15}{112} \cdot \frac{1}{70} = \frac{15}{7840}.
\]
Constraint E (the union of Constraints E1 and E2) is satisfied if and only if Constraint E2 is satisfied: assigning the four UEFA teams in Pot 2 to the four groups with at most one UEFA team implies that the remaining four non-UEFA teams in Pot 2 are all assigned to four of the six groups with a UEFA team, which is the necessary and sufficient condition for Constraint E1 in this case.

\item % b
\emph{One non-UEFA team in Pot 1 is assigned to one of the four groups without a UEFA team, and the other non-UEFA team in Pot 1 is assigned to one of the four groups with one UEFA team.}

The probability of this event is
\[
\frac{5}{8} \cdot \left( 2 \cdot \frac{4}{8} \cdot \frac{4}{7} \right) = \frac{5}{8} \cdot \frac{4}{7} = \frac{5}{14}.
\]
Then, three groups have two UEFA teams, and four groups have one UEFA team when Pot 1 is emptied.

Constraint E1 is satisfied if the four non-UEFA teams in Pot 2 are all assigned to four of the seven groups with a UEFA team, which has a probability of
\[
\frac{5}{14} \cdot \left( \frac{7}{8} \cdot \frac{6}{7} \cdot \frac{5}{6} \cdot \frac{4}{5} \right) = \frac{5}{14} \cdot \frac{1}{2} = \frac{5}{28}.
\]

Constraint E2 is satisfied if the four UEFA teams in Pot 2 are all assigned to four of the five groups with at most one UEFA team, which has a probability of
\[
\frac{5}{14} \cdot \left( \frac{5}{8} \cdot \frac{4}{7} \cdot \frac{3}{6} \cdot \frac{2}{5} \right) = \frac{5}{14} \cdot \frac{1}{14} = \frac{5}{196}.
\]

Constraint E (the union of Constraint E1 and E2) is violated if either Constraint E2 is violated, or the four UEFA teams in Pot 2 are all assigned to the four groups with one UEFA team. The latter event has a probability of
\[
\frac{5}{14} \cdot \left( \frac{4}{8} \cdot \frac{3}{7} \cdot \frac{2}{6} \cdot \frac{1}{5} \right) = \frac{5}{14} \cdot \frac{1}{70} = \frac{1}{196}.
\]
Thus, Constraint E is satisfied with a probability of 
\[
\frac{5}{196} - \frac{1}{196} = \frac{4}{196}.
\]

\item % c
\emph{Both non-UEFA teams in Pot 1 are assigned to two of the four groups with one UEFA team.}

The probability of this event is
\[
\frac{5}{8} \cdot \left( \frac{4}{8} \cdot \frac{3}{7} \right) = \frac{5}{8} \cdot \frac{3}{14} = \frac{15}{112}.
\]
Then, two groups have two UEFA teams, and six groups have one UEFA team when Pot 1 is emptied.

Constraint E1 is certainly satisfied.

Constraint E2 is satisfied if the four UEFA teams in Pot 2 are all assigned to four of the six groups with one UEFA team, which has a probability of
\[
\frac{15}{112} \cdot \left( \frac{6}{8} \cdot \frac{5}{7} \cdot \frac{4}{6} \cdot \frac{3}{5} \right) = \frac{15}{112} \cdot \frac{3}{14} = \frac{45}{1568}.
\]
Constraint E (the union of Constraints E1 and E2) is satisfied if and only if Constraint E2 is satisfied.
\end{enumerate}
\end{itemize}

To summarise, Constraints E1, E2 and E are satisfied with the following probabilities in the 2018 FIFA World Cup draw:
\begin{itemize}%[label=$\bullet$]
\item
Constraint E1:
\[
\frac{45}{1568} + \frac{15}{224} + \frac{15}{448} + \frac{15}{224} + \frac{3}{112} + \frac{45}{1568} + \frac{5}{28} + \frac{15}{112} = \frac{1769}{3136} \approx 56.41\%;
\]
\item
Constraint E2:
\[
0 + 0 + \frac{15}{3136} + 0 + \frac{9}{1568} + \frac{15}{7840} + \frac{5}{196} + \frac{45}{1568} = \frac{209}{3136} \approx 6.665\%;
\]
\item
Constraint E:
\[
0 + 0 + \frac{12}{3136} + 0 + \frac{9}{1568} + \frac{15}{7840} + \frac{4}{196} + \frac{45}{1568} = \frac{95}{1568} \approx 6.059\%.
\]
\end{itemize}

Our almost 162 million random unconstrained draws checked by the rejection sampler approximate well the last probability with 6.061\%; the error is less than 0.04\% or two thousandth percentage points.

\clearpage
\setcounter{table}{0}
\renewcommand{\thetable}{A.\arabic{table}}

\subsection*{Supplementary tables}
\addcontentsline{toc}{subsection}{Supplementary tables}

\begin{table}[ht!]
  \centering
  \caption{The non-uniformity of the 2018 FIFA World Cup draw}
  \label{Table_A1}
  
\begin{subtable}{\textwidth}
  \caption{Skip mechanism, the host (Russia) is pre-assigned to Group A (official draw procedure)}
  \label{Table_A1a}
\resizebox{\textwidth}{!}{
\begin{tiny}
\begin{tabularx}{1.43\textwidth}{r CCCC CCCC CCCC CCCC CCCC CCCC} \toprule
          & \rotatebox[origin=l]{90}{Spain} & \rotatebox[origin=l]{90}{Peru} & \rotatebox[origin=l]{90}{Switzerland} & \rotatebox[origin=l]{90}{England} & \rotatebox[origin=l]{90}{Colombia} & \rotatebox[origin=l]{90}{Mexico} & \rotatebox[origin=l]{90}{Uruguay} & \rotatebox[origin=l]{90}{Croatia} & \rotatebox[origin=l]{90}{Denmark} & \rotatebox[origin=l]{90}{Iceland} & \rotatebox[origin=l]{90}{Costa Rica} & \rotatebox[origin=l]{90}{Sweden} & \rotatebox[origin=l]{90}{Tunisia} & \rotatebox[origin=l]{90}{Egypt} & \rotatebox[origin=l]{90}{Senegal} & \rotatebox[origin=l]{90}{Iran} & \rotatebox[origin=l]{90}{Serbia} & \rotatebox[origin=l]{90}{Nigeria} & \rotatebox[origin=l]{90}{Australia} & \rotatebox[origin=l]{90}{Japan} & \rotatebox[origin=l]{90}{Morocco} & \rotatebox[origin=l]{90}{Panama} & \rotatebox[origin=l]{90}{South Korea} & \rotatebox[origin=l]{90}{Saudi Arabia} \\
\end{tabularx}
\end{tiny}
}
\end{subtable}

\vspace{0.25cm}
\begin{subtable}{\textwidth}
  \caption{Skip mechanism, the group with the host (Russia) is called Group A (our proposal)}
  \label{Table_A1b}
\resizebox{\textwidth}{!}{
\begin{threeparttable}
\begin{tiny}
    \begin{tabularx}{1.43\textwidth}{r CCCC CCCC CCCC CCCC CCCC CCCC} \toprule
          & \rotatebox[origin=l]{90}{Spain} & \rotatebox[origin=l]{90}{Peru} & \rotatebox[origin=l]{90}{Switzerland} & \rotatebox[origin=l]{90}{England} & \rotatebox[origin=l]{90}{Colombia} & \rotatebox[origin=l]{90}{Mexico} & \rotatebox[origin=l]{90}{Uruguay} & \rotatebox[origin=l]{90}{Croatia} & \rotatebox[origin=l]{90}{Denmark} & \rotatebox[origin=l]{90}{Iceland} & \rotatebox[origin=l]{90}{Costa Rica} & \rotatebox[origin=l]{90}{Sweden} & \rotatebox[origin=l]{90}{Tunisia} & \rotatebox[origin=l]{90}{Egypt} & \rotatebox[origin=l]{90}{Senegal} & \rotatebox[origin=l]{90}{Iran} & \rotatebox[origin=l]{90}{Serbia} & \rotatebox[origin=l]{90}{Nigeria} & \rotatebox[origin=l]{90}{Australia} & \rotatebox[origin=l]{90}{Japan} & \rotatebox[origin=l]{90}{Morocco} & \rotatebox[origin=l]{90}{Panama} & \rotatebox[origin=l]{90}{South Korea} & \rotatebox[origin=l]{90}{Saudi Arabia} \\
\end{tabularx}
\end{tiny}
    \begin{tablenotes} \footnotesize
\item
The probabilities are based on 10 million simulated draws.
\item
Blank cells represent pairs of teams that cannot play in the group stage.
\item
The numbers show percentages ($100 \times \Delta_{ij}$) rounded to one decimal place.
\item
\textcolor{ForestGreen}{Green} (\textcolor{red}{Red}) colour means that the Skip mechanism implies a higher (lower) probability than a uniform draw. For instance, the probability of assigning Russia and Nigeria to the same group is lower (higher) by about 2.1 (0.1) percentage points with the official draw mechanism (our proposal) compared to a uniform draw, when it is 12.5\%.
\item
Darker colour indicates a higher value.
    \end{tablenotes}
\end{threeparttable}
}
\end{subtable}
\end{table}

%\end{document}

%\input{Preambulum}

\begin{table}[ht!]
  \centering
  \caption{The non-uniformity of the 2022 FIFA World Cup draw}
  \label{Table_A2}
  
\begin{subtable}{\textwidth}
  \caption{Skip mechanism, the host (Qatar) is pre-assigned to Group A (official draw procedure)}
  \label{Table_A2a}
\resizebox{\textwidth}{!}{
\begin{tiny}
\begin{tabularx}{1.41\textwidth}{r CCCC CCCC CCCC CCCC CCCC CCCC} \toprule
           & \rotatebox[origin=l]{90}{Mexico} & \rotatebox[origin=l]{90}{Netherlands} & \rotatebox[origin=l]{90}{Denmark} & \rotatebox[origin=l]{90}{Germany} & \rotatebox[origin=l]{90}{Uruguay} & \rotatebox[origin=l]{90}{Switzerland} & \rotatebox[origin=l]{90}{United States} & \rotatebox[origin=l]{90}{Croatia} & \rotatebox[origin=l]{90}{Senegal} & \rotatebox[origin=l]{90}{Iran} & \rotatebox[origin=l]{90}{Japan} & \rotatebox[origin=l]{90}{Morocco} & \rotatebox[origin=l]{90}{Serbia} & \rotatebox[origin=l]{90}{Poland} & \rotatebox[origin=l]{90}{South Korea} & \rotatebox[origin=l]{90}{Tunisia} & \rotatebox[origin=l]{90}{Cameroon} & \rotatebox[origin=l]{90}{Canada} & \rotatebox[origin=l]{90}{Ecuador} & \rotatebox[origin=l]{90}{Saudi Arabia} & \rotatebox[origin=l]{90}{Ghana} & \rotatebox[origin=l]{90}{IPO1} & \rotatebox[origin=l]{90}{IPO2} & \rotatebox[origin=l]{90}{UEFA PO} \\ \bottomrule 
    Qatar & \cellcolor{ForestGreen!31.4787272727273} 3.5 & \cellcolor{red!10.5567272727273} 1.2 & \cellcolor{red!10.5214545454546} 1.2 & \cellcolor{red!10.52} 1.2 & \cellcolor{red!10.2854545454546} 1.1 & \cellcolor{red!10.7031818181818} 1.2 & \cellcolor{ForestGreen!31.4509090909091} 3.5 & \cellcolor{red!10.3428181818182} 1.1 & \cellcolor{ForestGreen!63.8406363636364} 7 &       &       & \cellcolor{ForestGreen!63.5729090909091} 7 & \cellcolor{red!95.3004545454546} 10.5 & \cellcolor{red!95.3680909090909} 10.5 &       & \cellcolor{ForestGreen!63.255} 7 & \cellcolor{red!38.9851818181818} 4.3 & \cellcolor{ForestGreen!9.31563636363635} 1 & \cellcolor{red!27.6128181818182} 3 &       & \cellcolor{red!38.9942727272727} 4.3 &       & \cellcolor{ForestGreen!9.05672727272728} 1 & \cellcolor{ForestGreen!87.2199090909091} 9.6 \\
\end{tabularx}
\end{tiny}
}
\end{subtable}

\vspace{0.25cm}
\begin{subtable}{\textwidth}
  \caption{Skip mechanism, the group with the host (Qatar) is called Group A (our proposal)}
  \label{Table_A2b}
\resizebox{\textwidth}{!}{
\begin{threeparttable}
\begin{tiny}
    \begin{tabularx}{1.41\textwidth}{r CCCC CCCC CCCC CCCC CCCC CCCC} \toprule
           & \rotatebox[origin=l]{90}{Mexico} & \rotatebox[origin=l]{90}{Netherlands} & \rotatebox[origin=l]{90}{Denmark} & \rotatebox[origin=l]{90}{Germany} & \rotatebox[origin=l]{90}{Uruguay} & \rotatebox[origin=l]{90}{Switzerland} & \rotatebox[origin=l]{90}{United States} & \rotatebox[origin=l]{90}{Croatia} & \rotatebox[origin=l]{90}{Senegal} & \rotatebox[origin=l]{90}{Iran} & \rotatebox[origin=l]{90}{Japan} & \rotatebox[origin=l]{90}{Morocco} & \rotatebox[origin=l]{90}{Serbia} & \rotatebox[origin=l]{90}{Poland} & \rotatebox[origin=l]{90}{South Korea} & \rotatebox[origin=l]{90}{Tunisia} & \rotatebox[origin=l]{90}{Cameroon} & \rotatebox[origin=l]{90}{Canada} & \rotatebox[origin=l]{90}{Ecuador} & \rotatebox[origin=l]{90}{Saudi Arabia} & \rotatebox[origin=l]{90}{Ghana} & \rotatebox[origin=l]{90}{IPO1} & \rotatebox[origin=l]{90}{IPO2} & \rotatebox[origin=l]{90}{UEFA PO} \\ \bottomrule
    Qatar & \cellcolor{ForestGreen!26.2778181818182} 2.9 & \cellcolor{red!15.8115454545455} 1.7 & \cellcolor{red!16.0450909090909} 1.8 & \cellcolor{red!16.1093636363636} 1.8 & \cellcolor{ForestGreen!27.5739090909091} 3 & \cellcolor{red!16.1617272727273} 1.8 & \cellcolor{ForestGreen!26.3032727272727} 2.9 & \cellcolor{red!16.0272727272727} 1.8 & \cellcolor{ForestGreen!22.7908181818182} 2.5 &       &       & \cellcolor{ForestGreen!22.6740909090909} 2.5 & \cellcolor{red!33.9007272727273} 3.7 & \cellcolor{red!34.1079090909091} 3.8 &       & \cellcolor{ForestGreen!22.5437272727273} 2.5 & \cellcolor{red!9.70981818181817} 1.1 & \cellcolor{red!10.4170909090909} 1.1 & \cellcolor{red!11.2981818181818} 1.2 &       & \cellcolor{red!9.79081818181819} 1.1 &       & \cellcolor{red!10.4659090909091} 1.2 & \cellcolor{ForestGreen!51.6818181818182} 5.7 \\
\end{tabularx}
\end{tiny}
    \begin{tablenotes} \footnotesize
\item
The probabilities are based on 10 million simulated draws.
\item
X represents a pair of teams that cannot play in the group stage.
\item
The numbers show percentages ($100 \times \Delta_{ij}$) rounded to one decimal place.
\item
\textcolor{ForestGreen}{Green} (\textcolor{red}{Red}) colour means that the Skip mechanism implies a higher (lower) probability than a uniform draw. For instance, the probability of assigning Qatar and Uruguay to the same group is lower (higher) by about 1.1 (3.0) percentage points with the official draw mechanism (our proposal) compared to a uniform draw, when it is 13.6\%.
\item
Darker colour indicates a higher value.
    \end{tablenotes}
\end{threeparttable}
}
\end{subtable}
\end{table}

%\end{document}

\end{document}